\documentclass[psamsfonts]{amsart}

\usepackage{amssymb,amsfonts}
\usepackage{soul}
\usepackage{amsmath}
\usepackage[all,arc]{xy}
\usepackage{enumerate}
\usepackage{mathrsfs}
\usepackage{color}
\newtheorem{thm}{Theorem}[section]
\newtheorem{cor}[thm]{Corollary}
\newtheorem{prop}[thm]{Proposition}
\newtheorem{lem}[thm]{Lemma}

\newtheorem{quest}[thm]{Question}
\newtheorem{claim}[thm]{Claim}
\newtheorem{defn}[thm]{Definition}
\newtheorem*{quest*}{Question}

\theoremstyle{definition}

\theoremstyle{remark}
\newtheorem{rem}[thm]{Remark}

\newcommand{\Z}{\mathbb{Z}}

\newcommand{\x}{\mathcal{X}} 
\newcommand{\xp}{E} 
\newcommand{\xu}{E_n} 
\newcommand{\xut}{\widetilde{\xu}} 
 
\newcommand{\xf}{\xn_\text{\rm flex}}
 
\newcommand{\xn}{\widetilde{\x}} 
 
\newcommand{\pf}{p^\mathrm{flex}} 

\newcommand{\pgl}{\mathrm{PGL}_3(\CC)} 
\newcommand{\ok}{\overline{K}} 

\newcommand{\pix}{\pi_1(\x,F)}

\newcommand{\ZZ}{\mathbb{Z}} 
\newcommand{\HH}{\mathbb{H}} 
\newcommand{\RR}{\mathbb{R}} 
\newcommand{\CC}{\mathbb{C}} 
\newcommand{\QQ}{\mathbb{Q}} 
\newcommand{\PP}{\mathbf{P}} 
\newcommand{\slz}{\mathrm{SL}_2(\ZZ)} 
\newcommand{\slr}{\mathrm{SL}_2(\RR)} 
 
\newcommand{\sflex}{s_\text{flex}}

\newcommand{\slc}{\mathrm{SL}_3(\CC)} 
\newcommand{\gm}{\Gamma}
\newcommand{\gn}{{\gm}_n}
\newcommand{\gt}{\Lambda}
\newcommand{\gf}{{\gm_\mathrm{flex}}} 
\newcommand{\uconf}{\mathrm{UConf}}

\newcommand{\TT}{\mathcal{T}_{1,1}} 
\newcommand{\CT}{\mathcal{C}_{1,1}}

\usepackage{graphicx}

\newcommand{\xm}{\widetilde{\x}_m}

\title[Choosing points on cubic plane curves]{Choosing points on cubic plane curves: rigidity and flexibility}

\author{Ishan Banerjee}
\address{Ishan Banerjee, Department of Mathematics, University of Chicago, 5734 S University Ave, Chicago, IL, USA.}
\email{ishan@math.uchicago.edu}

\author{Weiyan Chen}
\address{Weiyan Chen, Yau Mathematical Sciences Center, Ning Zhai, Tsinghua University, Hai Dian District, Beijing, China.}
\email{chwy@tsinghua.edu.cn}

\subjclass[2010]{Primary 55R10; Secondary 55R80, 14H50.}

\begin{document}

\begin{abstract}
    Every smooth cubic plane curve has 9 flex points and 27 sextatic points. We study the following question of Farb: Is it true that the known algebraic structures give all the possible ways to continuously choose $n$ distinct points on every smooth cubic plane curve, for each given positive integer $n$? We give an affirmative answer to the question when $n=9$ and $18$ (the smallest open cases), and a negative answer for infinitely many $n$'s.
\end{abstract}

\maketitle

\section{Introduction}

It is a classical topic to study certain structures of special points on complex smooth cubic plane curves, for example, the 9 flex points (first attributed to Maclaurin; see Introduction of \cite{AD} for a brief history) and the 27 sextatic points (studied by Cayley \cite{Cayley}). 
Inspired by these algebro-geometric constructions, Benson Farb asked the following  question:
\begin{quest}\label{question}
Is it true that the known algebraic structures give all the possible ways to continuously choose $n$ distinct points on every smooth cubic plane curve as the curve varies in families, for any given positive integer $n$? 
\end{quest}


Let us make precise what it means to ``continuously choose  $n$  points" in Question \ref{question}. Let $\x$ denote the space of all smooth cubic plane curves:
$$\x = \{ F(x,y,z)\ |\ F \textrm{ is a nonsingular homogeneous polynomial of degree 3}  \} / \CC^\times.$$
Let $\pi:\xp\to\x$ denote the universal cubic plane curve, a fiber bundle whose  fiber over $F\in\x$ is the smooth cubic plane  curve $C_F$ defined by the equation $F=0$. A \emph{multisection} of $\pi$ of degree $n$ is a continuous section of an associated bundle whose fiber over $F\in\x$ is $\uconf_n(C_F)$, the space of $n$ distinct unordered points in $C_F$:
$$\xymatrix{
\uconf_n(C_F) \ar[r] & \xu \ar[r] & \x. \ar@{-->}@/^-10pt/[l]\\
}$$
To choose  $n$ distinct points continuously on every cubic curve is to give a multisection of $\pi$ of degree $n$. Constructions from algebraic geometry give the following:
\begin{thm}[Maclaurin, Cayley, Gattazzo; see Section \ref{alg cons} for exact attributions]
\label{old AG}
The universal cubic plane curve $\pi$ admits a multisection of degree $n$ when $$n=9\sum_{m\in I} J_2(m), \ \ \ \ \ \  \ \ \ \ \text{ where } J_2(m)=m^2\prod_{p|m,\text{ prime}}(1-p^{-2})$$ 
and $I$ is a set of positive integers, for example, when $n= 9, 27, 36, 72, 81, 99, 108...$
\end{thm}
How are these multisections constructed? Every smooth cubic plane curve $C_F$ defined by $F=0$ is a Riemann surface of genus 1. We can pick all the $3m$-torsion points on $C_F$ as an elliptic curve with any flex point chosen as identity. This construction is independent of our choice of flex points since they differ by 3-torsion. We can also take the disjoint union, among $m\in I$, of such constructions. All multisections in Theorem \ref{old AG} are constructed in this way. We will call them \emph{multisections from torsion construction}. See Section \ref{alg cons} for more detail. In this setting, Question \ref{question} asks: Is it true that multisections from torsion construction are the only {continuous} multisections of $\pi$ up to homotopy? 

The first step towards answering this question was made by the second author in \cite{WC}:
\begin{thm}[Chen, Theorem 2 in \cite{WC}]\label{old thm}
If $9\nmid n$, then the universal cubic plane curve $\pi$ does not have any continuous multisection of degree $n$.
\end{thm}
This theorem answered Question \ref{question} affirmatively when $9\nmid n$. In the present paper, we continue to study Question \ref{question} when $n$ is a multiple of 9. We proved:
\begin{thm}\label{9}
Any multisection of $\pi$ of degree $9$ must be homotopic to the multisection given by choosing the $9$ flex points on each cubic curve. 
\end{thm}

\begin{rem}\label{why 9}
The flex points on cubic curves have a delicate structure that is well-studied. They form the \emph{Hesse configuration} of 9 points and 12 lines (Hesse \cite{Hesse} \cite{Hesse2}, 1844) whose monodromy group is a finite group of order 216 called the \emph{Hesse group} (Jordan \cite{Jordan}, 1877). Theorem \ref{9} tells us that this delicate structure is also unique: It is the only
continuous choice of 9 distinct points (up to homotopy) on all smooth cubic curves. We thus obtain an affirmative answer to Question \ref{question} for $n=9$, the first open case after Theorem \ref{old thm}.
\end{rem}

When $n=18$, we have the following:
\begin{thm}\label{18}
The universal cubic plane curve $\pi$ has no multisection of degree $18$.
\end{thm}
\begin{rem}
\label{why 18}
Notice that 18 is not among those $n$'s that appeared in Theorem \ref{old AG}. Algebraic geometers have not discovered any structure of 18 marked points on cubic curves. There is, as Theorem \ref{18} shows, a topological reason for that. 

\end{rem}

Theorem \ref{old thm}, \ref{9} and \ref{18} above all give positive answers to Question \ref{question} for various  $n$'s. One might wonder whether the answer to the question is ``yes" for every integer $n$. This, as we will prove, is not true. We will produce an infinite list of degree $n$'s for which Question \ref{question} has a negative answer. The idea is to deform an existing multisection from torsion construction and obtain new multisections of various different degrees. 
\begin{thm}\label{connected new multi}
For any $m\ge 4$, the universal cubic plane curve $\pi$ has a connected multisection of degree $n=18J_2(m)$ that is not homotopic to any multisection from torsion construction (\emph{i.e.} those multisections in Theorem \ref{old AG}).
\end{thm}
Here a multisection is \emph{connected} if the space $\xn$ of smooth cubic curves equipped with a point in the multisection is connected ($\xn$ is a cover of $\x$). See Section \ref{section mult vs vir} for more details. Theorem \ref{connected new multi} gives negative answers to Question \ref{question} for all those $n$'s as stated, such as $n=$ 216, 432, 864, 1296, 2160... The reason for the hypothesis $m\ge4$ is to ensure that  the principal congruence subgroup $\Gamma_1(m)$ in $\slz$ is free, a condition needed in the proof.


Finally, if we allow our multisection to be disconnected, then we can construct many new multisections in the following theorem.
\begin{thm}
\label{topo constr}
The universal cubic plane curve $\pi$ has a multisection of degree $n$ that is not homotopic to any multisection from torsion construction (\emph{i.e.} those in Theorem \ref{old AG}) when $n$ can be expressed as 
$$n=9\sum_{m\in I} k_m J_2(m)$$
where $I$ is a finite set of positive integers and each $k_m$ is a positive integer such that $k_m=1$ for every $m\le3$ and that $k_m>1$ for some $m\ge 4.$ 
\end{thm}
This theorem gives negative answers to Question \ref{question} for all those $n$'s as stated, such as $n=$ 225, 243, 253, 288, 297, 315, 324, 333, ... 

Our theorems above can be roughly summarized in one slogan: multisections of low degrees are rigid (Theorem \ref{9} and \ref{18}), while multisections of high degrees are flexible as they can be deformed to new ones (Theorem \ref{connected new multi} and \ref{topo constr}).

\begin{rem}[\textbf{Strategy and difficulties}]
Our strategy for proofs is to relate the problem about the universal cubic plane curve $\pi$ to the same problem about the universal torus bundle over the classifying space $\mathrm{B}\slz$ where $\slz$ is the mapping class group of the torus. However, there are two main difficulties. First, the bundle $\pi$ is not a pullback of the universal torus bundle over $\mathrm{B}\slz$ because $\pi$ does not have any continuous section (Theorem \ref{old thm} for $n=1$). Second, even if $\pi$ were a pullback from a bundle over $\mathrm{B}\slz$, nonzero obstruction classes on $\mathrm{B}\slz$ might pullback to zero classes on $\x$, which is not useful to obstruct multisections from $\x$. We will overcome the first difficulty by studying the fundamental groups of covers of $\x$ associated to multisections, and overcome the second difficulty by cohomological computations using spectral sequences. 
\end{rem}

\begin{rem}[\textbf{Algebraic multisections}]
One could also ask Question \ref{question} for \emph{algebraic} multisections instead of continuous multisections. For this variant question, it turns out that any \emph{algebraic} multisection must be a multisection from torsion construction. Curt McMullen  told us this fact and sketched a proof to us, after we announced an earlier draft of the present paper. Since McMullen's observation answers an interesting variant of Question \ref{question} which we often got asked by the audience when giving a talk on this subject, we think it is worthwhile to record it in an Appendix to this paper, with McMullen's gratefully acknowledged authorization.

All results in the rest of the paper are about continuous multisections and therefore are logically independent of the  result in the Appendix. However, Theorem \ref{algebraic is torsion} in the Appendix implies the following strengthening of Theorem \ref{connected new multi} and Theorem \ref{topo constr}: The universal cubic plane curve $\pi$ has a \emph{continuous} multisection that is not homotopic to any \emph{algebraic} multisections.
\end{rem}

The paper proceeds as following: In Section 2, we compare the notions of multisections and virtual sections. In Section 3 and 4, we prove several preliminary lemmas about the fundamental groups and the cohomology groups of covers of $\x$ associated to multisections. In Section 5, we study obstructions to multisections and prove Theorem \ref{9} and \ref{18}. In Section 6, we discuss algebraic constructions of multisections that are previously known. In Section 7, we give new topological constructions of multisections and prove Theorem \ref{connected new multi} and \ref{topo constr}. In the Appendix, we record McMullen's argument showing that every algebraic multisection must be from torsion construction.



\section*{Acknowledgment}
We are grateful to Benson Farb for sharing his question. We thank Igor Dolgachev, Jordan Ellenberg, and Tom Goodwillie for helpful conversations, and thank Xiaoxiao Huang for her help in drawing those figures in the paper. We especially thank Curt McMullen for explaining to us a related result about algebraic multisections and allowing us to record his observation in the Appendix. The second author is grateful to the School of Mathematics at the University of Minnesota Twin Cities for its hospitality, as most part of this joint work was done when he was a postdoc there. 

\section*{Notations}
Here is a list of frequently used notations with the sections where they first appear.
\begin{enumerate}
    \item $\x$ denotes the space of smooth cubic plane curves. (Introduction)
    \item $\xp\xrightarrow[]{\pi}\x$ denotes the universal cubic curve, which is a fiber bundle with base $\x$ whose fiber over $F\in\x$ is the cubic curve $C_F.$ (Introduction)
    \item We generally use $p:\xn\to\x$ to denote the cover induced by a multisection. (Section \ref{section mult vs vir})
    \item In particular, we use $p_m:\xm\to\x$ to denote the cover induced by the multisection given by choosing points of type $3m$ for $m\in\ZZ_{>0}$. (Section \ref{alg cons})
    \item $\pf:\xf\to \x$ denotes the flex cover of $\x$. (Section \ref{section mult vs vir})
    \item $K$ denotes the kernel of the monodromy representation $\pi_1(\x)\to \slz$ of the bundle $\pi$. Its center is $Z(K).$ 
    (Section \ref{fund group})
    \item $\gm:=\pi_1(\x)/Z(K)$, while $\gn:=\pi_1(\xn)/Z(K)$. (Section \ref{fund group})
    
\end{enumerate}

\section{Multisections and virtual sections} 
\label{section mult vs vir}

In this section, we will compare \emph{multisections} and \emph{virtual sections} as two notions of ``choosing $n$ points continuously", following the terminologies of L. Chen and Salter \cite{CS}  in a different context. While we phrased our main theorems in terms of multisections in the Introduction, our proofs later will focus on virtual sections because they have more structures to work with. Therefore, it is important to understand how these two notions are related. 

As before, the universal cubic plane curve is the fiber bundle $\pi: \xp\to\x$ whose fiber over $F\in\x$ is the curve $C_F$ defined by the equation $F=0$. 

\begin{defn} [\bf{Virtual section}]\label{def virtual}
A {virtual section} of $\pi$ of degree $n$ is a triple $(\xn,p,s)$ where  $p: \xn \to \x$ is a (possibly disconnected) $n$-sheeted cover  and  $s$ is a continuous map making the following diagram commute:
$$
\xymatrix{
& \xp \ar[d]^\pi\\
\xn \ar[r]^p \ar@{-->}[ru]^s & \x\\
}
$$


Two virtual sections $(\xn,p, s)$ and $(\xn',p',s')$ are {homotopic as virtual sections} if there is an isomorphism $f:\xn\to\xn'$ of covers of $\x$ and a homotopy $s_t:\xn\to \xp$ such that $s_t$ is a virtual section for each $t\in[0,1]$ and $s_0=s'\circ f$ and $s_1=s$. 

\end{defn}

For any topological space $Y$, define 
$$\uconf_n(Y):=\{(y_1,\cdots,y_n)\in Y^n \ | \ y_i\ne y_j \text{ when $i\ne j$}\}/S_n$$
where the symmetric group $S_n$ acts on the $n$-tuples by permuting the ordering. $\uconf_n(Y)$ is the configuration space of $n$ distinct unordered points in $Y$. Each point in $\uconf_n(Y)$ is an $n$-element subset $\{y_1,\cdots,y_n\}$ of $Y$. 


\begin{defn}[\bf{Multisection}]\label{multi def}
A {multisection} of degree $n$ of the universal cubic curve $\pi$ is a continuous section $\sigma$ of the fiber bundle $\uconf_n(C_F)\to \xu\to \x$ whose fiber over $F\in\x$ is $\uconf_n(C_F)$:
$$\xymatrix{
\uconf_n(C_F) \ar[r] & \xu \ar[d] \\
& \x \ar@{-->}@/^-10pt/[u]_\sigma \\
}$$
Two multisections $\sigma_1$ and $\sigma_2$ of the same degree are {homotopic as multisections} if they are homotopic as sections of $\xu\to\x.$
\end{defn}
Although one can define virtual sections or multisections of any fiber bundle,  
in the present paper we will only discuss multisections or virtual sections of the universal cubic curve $\pi$. To be brief, we might not explicitly mention $\pi$ each time we mention virtual sections or multisections.

Next we study the relationship between virtual sections and multisections.

\begin{prop}
\label{multi gives virtual}
Every multisection induces a virtual section of the same degree.
\end{prop}
\begin{proof}
Suppose $\sigma:\x\to \xu$ is a multisection of degree $n$. Define the following decorated configuration space of $n$ distinct unordered points with one point marked on $C_F$
\begin{equation}
    \label{conf_n,1}
    \uconf_{n,1}(C_F):=\{(S,x)\in \uconf_n(C_F)\times C_F : x\in S\}.
\end{equation}
The projection map $\uconf_{n,1}(C_F)\to\uconf_n(C_F)$ onto the first coordinate is a covering map of degree $n$. Applying this construction fiberwise, we obtain the following fiber bundle 
$$\uconf_{n,1}(C_F)\to \xut\to \x.$$
The total space $\xut$ is a cover of $\xu$ of degree $n$. Now, define  $\xn:=\sigma^* \xut$, the pullback of the cover $\xut\to\xp$ via the multisection $\sigma$. 
More concretely, the space $\xn$ consists of pairs $(F,x)$ where $F\in \x$ and $x$ is a point in the set $\sigma(F)$ of $n$ distinct points on $C_F$. The covering map $p:\xn\to \x$ is via $(F,x)\mapsto F$. Moreover, since $x\in C_F$, the cover $\xn$ is naturally a subspace of $\xp$. The inclusion $s:\xn\hookrightarrow\xp$ gives the virtual section. 
\end{proof}

In the setting of the proof above, we will  call $\xn$ the \emph{cover of $\x$ associated to the multisection $\sigma$}. We  say a multisection is \emph{connected} if its associated cover is connected.

The following proposition gives a partial converse to Proposition \ref{multi gives virtual}. 

\begin{prop}
\label{multi is virtual}
A virtual section $(\xn,p,s)$ induces a multisection of the same degree if $s$ is injective. 


\end{prop}
\begin{proof}
For each $F\in \x$, suppose $p^{-1}(F)=\{F_1,\cdots, F_n\}$. Then $s(F_1),\cdots, s(F_n)$ are points on $C_F$ by Definition \ref{def virtual}. Since $s$ is injective, $s(F_1),\cdots,s(F_n)$ are all distinct. Hence, the map $\sigma: F\mapsto \{s(F_1),\cdots, s(F_n)\}$ defines a multisection of degree $n$.
\end{proof}
An injective virtual section $s$ identifies $\xn$ as a subspace of $\xp$. In this case, each point in $\xn$ can be uniquely represented as a pair $(F,x)$ where $F\in\x$ and $x\in C_F$.

It is straightforward to check that the construction of $s$ from $\sigma$ in Proposition \ref{multi gives virtual} and the construction of $\sigma$ from $s$ in Proposition \ref{multi is virtual} are inverses of each other. Hence, the two propositions together give the following corollary:
\begin{cor}
\label{multi is inj virtual}
Multisections are in bijection with injective virtual sections.
\end{cor}
Our readers should feel free to think  multisections and injective virtual sections as the same. However, we will use $\sigma$ (or other Greek letters) to denote multisections and use $s$ to denote virtual sections because, technically speaking, they are functions from and to different spaces.

Next, we study the relation between the notions of homotopy through multisections and homotopy through virtual sections. 

\begin{prop}\label{mult vs virtual homotopy1}
Suppose that $\sigma_i$ is a multisection of $\pi$ that induces the virtual section $s_i$ for each $i=0,1$. If $\sigma_0$ and $\sigma_1$ are homotopic as multisections, then $s_0$ and $s_1$ are homotopic as virtual sections.  
\end{prop}
\begin{proof}
The proposition should be  intuitively clear: By Corollary \ref{multi is inj virtual}, a homotopy through multisections is a homotopy through injective virtual sections. To be rigorous, we include a careful argument below.

Suppose that $\sigma_t$ for $t\in[0,1]$ is a homotopy between $\sigma_0$ and $\sigma_1$. As in the proof of Proposition \ref{multi gives virtual}, for each $t\in[0,1]$, define 
$$\xn_t:=\Big\{(F,x): F\in\x \text{ and } x\in \sigma_t(F)\Big\}.$$
The map $p_t:\xn_t\to\x$ given by $(F,x)\mapsto F$ is an $n$-sheeted cover at any $t$. By our definition, the space $\xn_t$  is already a subspace of $\xp$ and hence the inclusion defines a virtual section $s_t$. We claim that the family $(\xn_t,p_t,s_t)$ gives a homotopy between the virtual sections $(\xn_0,p_0,s_0)$ and $(\xn_1,p_1,s_1)$. 

We first construct an isomorphism of covers $f:\xn_0\to\xn_1$ as needed in Definition \ref{def virtual}.  Take any $(F,x)\in \xn_0$. 
Consider the function from $[0,1]$ to $\uconf_n(C_F)$ given by $t\mapsto \sigma_t(F)$. There is a unique lift of this map to the covering space $\uconf_{n,1}(C_F)$ defined in (\ref{conf_n,1}) when we choose the appropriate base points as below:
$$
\xymatrix{
& &  \bigg(\uconf_{n,1}C_F,\ \big(\sigma_0(F),x\big)\bigg)\ar[d]\\
\bigg([0,1],\ 0\bigg) \ar[rr]_{t\mapsto {\sigma}_t(F)} \ar@{-->}[rru]^{t\mapsto \widetilde{\sigma_t}(F)} & &\bigg(\uconf_nC_F,\  \sigma_0(F)\bigg) \\
}$$
Intuitively, the function $\widetilde{\sigma_t}(F)$ traces what the homotopy $\sigma_t$ does to the marked point $x$. In the end, $\widetilde{\sigma_1}(F)$ picks a marked point $y$ in $\sigma_1(F)$. The function 
\begin{align*}
    f:\xn_0&\to\xn_1\\
    (F,x)&\mapsto (F,y)
\end{align*}
gives an isomorphism of covers. $f$ is continuous because the homotopy $\sigma_t$ is continuous.

More generally, for any $t\in[0,1]$, we have an isomorphism of covers $f_t:\xn_0\to\xn_t$ by the same argument above with the end point $1$ replaced by $t$. Now the map $s_t\circ f_t: \xn_0\to\xp$ gives the desired homotopy of virtual sections as in Definition \ref{def virtual}.

\end{proof}

It turns out that the converse of Proposition \ref{mult vs virtual homotopy1} is also true. This converse is nontrivial because we need to show that a homotopy through virtual sections can give a homotopy through injective ones. The proof, though nontrivial, is not difficult and uses a standard transversality argument.

\begin{prop}\label{mult vs virtual homotopy2}
Suppose that $\sigma_i$ is a multisection of $\pi$ that induces the virtual section $s_i$ for each $i=0,1$. If $s_0$ and $s_1$ are homotopic as virtual sections, then $\sigma_0$ and $\sigma_1$ are homotopic as multisections.  
\end{prop}
\begin{proof}
Suppose that $s_t$ for $t\in[0,1]$ is a homotopy of virtual sections as in Definition \ref{def virtual}. Our goal is to use a transversality argument to homotope $s_t$ to another homotopy $h_t$ which is an injective virtual section for any $t\in[0,1]$, and hence a multisection by Corollary \ref{multi is inj virtual}. 

As before, let $\xn$ denote the cover associated to the multisection $\sigma_0$. Each point in $\xn$ can be uniquely represented as a pair $(F,x)$ where $F\in\x$ and $x\in C_F$. Define another cover $Y$ of $\x$ as 
$$Y:=\{(F,x_1,\cdots, x_n)\ |\ (F,x_i)\in \xn \text{ for all $i$ and }  x_i\ne x_j \text{ for all $i\ne j$}\}.$$
$Y$ is nothing but $\xn$ equipped with an ordering on its fiber over $\x$. More formally, $Y$ is the normalization of $\xn$ so that $\pi_1(Y)$ is the kernel of the monodromy representation $\pi_1(\xn)\to S_n$ of the cover $\xn\to\x$. Define 
$$E^{(n)}:=\{(F,x_1,\cdots, x_n)\ |\  x_i\in C_F \text{ for all $i$}\}.$$
The homotopy $s_t$ of virtual sections gives the following continuous map
\begin{align*}
    S:[0,1]\times Y &\longrightarrow E^{(n)}\\
    (t, F,x_1,\cdots,x_n)&\longmapsto (F, s_t(F,x_1), \cdots, s_t(F,x_n)).
\end{align*}
In the equation above, we've slightly abused the notation: $s_t(F,x)$ as an element in $\xp$ should be of the form $(F,y)$. But for simplicity we will drop the repeated $F$ from the notation and just use $s_t(F,x)$ to denote the point $y$. We can assume that $S$ is a $C^\infty$ function after a small homotopy relative to the boundary $\{0,1\}\times Y$ of its domain.

The symmetric group $S_n$ acts on both $Y$ and $E^{(n)}$ by permuting the ordering of the $x_i$'s. Let $\mathcal{F}$ denote the space of $C^\infty$ functions $H:[0,1]\times Y\to E^{(n)}$ such that 
\begin{enumerate}
    \item $H$ is $S_n$-equivariant
    \item $H=S$ on the boundary $\{0,1\}\times Y$
    \item the following diagram commutes where the two vertical maps denote natural projections onto $F\in\x$:
$$
\xymatrix{
[0,1]\times Y \ar[rr]^H \ar[rd] & & E^{(n)} \ar[ld]\\
& \x & \\
}
$$
\end{enumerate}
Note that  $S$ is an element in $\mathcal{F}$. 

\begin{claim}
\label{transverse}
We can homotope $S$ in the function space $\mathcal{F}$ to another function $H:[0,1]\times Y\to E^{(n)}$ such that $H$ is transverse to the following submanifold $\Delta E_{ij}$ of $E^{(n)}$ for every pair $i\ne j \in\{1,\cdots,n\}$:
$$\Delta E_{ij}:=\{(F, x_1,\cdots,x_n)\ |\ x_i= x_j\}\  \subset\   E^{(n)}.$$
\end{claim}
We can pick such an $H$ because transverse maps are generic.  However, since $\mathcal{F}$ has infinite dimension, we will need to use an infinite-dimensional version of the transversality theorem due to Smale \cite{Smale}. We will postpone more details to the end. Let's now assume the claim and finish the proof of the proposition.
\begin{figure}[h]
    \centering
    \includegraphics[scale=.36]{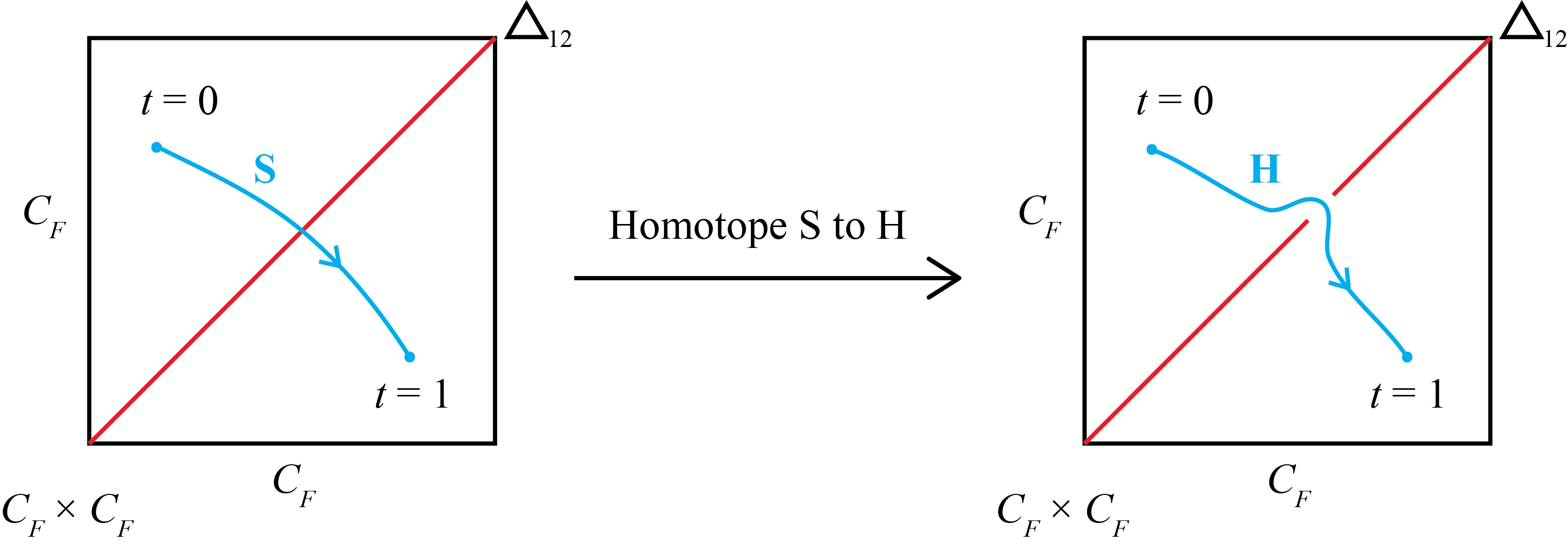}
     \caption{The figure illustrates the images of $S$ and  $H$ on a single fiber  of the bundle $E^{(n)}\to \x$ over a point $F$ when $n=2$. In this case, the fiber is $C_F\times C_F$. Since $C_F$ is 2-dimensional, $S$ can always be homotoped to an $H$ whose image has empty intersection with the diagonal $\Delta_{12}\subset C_F\times C_F$.}
     \label{transverse fig}
\end{figure}

The condition that $H$ is transverse to $\Delta E_{ij}$ implies that the image of $H$ has empty intersection with $\Delta E_{ij}$ by a simple dimension count: Since $H$ commutes with the bundle projections as in (3) above, we should focus on how $H$ maps fibers to fibers when counting dimensions. The fiber of the bundle projection $[0,1]\times Y\to \x$ is  1-dimensional. The fiber of the bundle projection $\Delta E_{ij}\to \x$ is of codimension 2 in the fiber of the bundle projection $E^{(n)}\to\x$. Since $H$ is transverse to $\Delta E_{ij}$, the image of $H$ must have empty intersection with $\Delta E_{ij}$. 
See Figure \ref{transverse fig} for a picture of the homotopy. Therefore, the map $H\in\mathcal{F}$ has the following representation:
\begin{align*}
    H:[0,1]\times Y &\longrightarrow E^{(n)}\setminus \bigcup_{i\ne j} \Delta E_{ij} \\
    (t, F,x_1,\cdots,x_n)&\longmapsto (F, h_t^1, \cdots, h_t^n).
\end{align*}
where each coordinate  $h^i_t=h^i_t(F,x_1,\cdots,x_n)$ is a point on $C_F$. Now, for each $t\in[0,1]$ and every $F\in\x$, define the set 
$$\sigma_t(F):=\{h^i_t(F,x_1,\cdots,x_n)\ |\ i=1,\cdots,n,\ \ (F,x_1,\cdots,x_n)\in Y\}.$$
A priori, the definition of $\sigma_t(F)$ depends not only on $F\in\x$ but also on an ordering of the points $(x_1,\cdots, x_n)$. However, since $H$ is $S_n$-equivariant by (1), reordering the $x_i$'s amounts to reordering the $h^i_t$'s without changing the set $\sigma_t(F)$. Hence, the set $\sigma_t(F)$ only depends on $F$. Moreover, since all the $h_t^i$'s are pairwise distinct,  $\sigma_t(F)$ is always an $n$-element subset of $C_F$. In other words, $\sigma_t: \x\to\xu$ is a multisection for all $t\in[0,1]$. The condition (2) that $H=S$ on $\{0,1\}\times Y$ implies that $\sigma_t$ is a homotopy of the multisections $\sigma_0$ and $\sigma_1$ that we started with. 

Finally, let us go back to Claim \ref{transverse} and briefly explain why we can choose such a transverse $H$. First of all, the evaluation $ev: \mathcal{F}\times ([0,1]\times Y)\to E^{(n)}$ given by $(H,y)\mapsto H(y)$ is a submersion. This is because given $n$ tangent vectors of $C_F$ on $n$ distinct points, a given $H$ in $\mathcal{F}$ can always be deformed along those $n$ tangent directions at the $n$ points while staying in $\mathcal{F}$. Hence, the preimage $ev^{-1}(\Delta E_{ij})$ is an (infinite-dimensional) submanifold of $\mathcal{F}\times ([0,1]\times Y)$. Consider the projection map $proj: ev^{-1}(\Delta E_{ij})\to \mathcal{F}$. By an infinite dimensional version of Sard's theorem proven by Smale (Theorem 1.3 in \cite{Smale}), the subset 
$$R_{ij}:=\{\text{regular values of $proj: ev^{-1}(\Delta E_{ij})\to \mathcal{F}$}\}$$
is residual (\emph{i.e.} has meager complement) in $\mathcal{F}$. It is straightforward to check that those regular values $H\in R_{ij}$ are precisely functions in $\mathcal{F}$ that are transverse to $\Delta E_{ij}$. The finite intersection $\bigcap_{i\ne j}R_{ij}$ is also residual in $\mathcal{F}$. Thus, we can always homotope $S$ while staying in $\mathcal{F}$ to an element $H\in \bigcap_{i\ne j}R_{ij}$ as in Claim \ref{transverse}.
\end{proof}

\begin{rem}[\textbf{Homotopy classes of multisections $\ne$ homotopy classes of virtual sections}]
Proposition \ref{mult vs virtual homotopy1} and \ref{mult vs virtual homotopy2} do not imply that there is a natural bijection between homotopy classes of multisections and homotopy classes of virtual sections. In fact, there does exist a virtual section that is not homotopic to any injective virtual section, as a consequence of one of our main theorems. For example, consider 
\begin{equation}
\label{def xf}
    \xf:=\{(F,x)\ |\ x \text{ is a flex point on $C_F$}\}.
\end{equation}
The natural inclusion $\sflex:\xf\to\xp$ gives a virtual section of degree 9. Now let $\xn$ be the disjoint union of two copies of $\xf$. We have a virtual section $s:\xn\to\xp$ by mapping each component of $\xn$ to $\xp$ via the same $\sflex$. The virtual section $s$ has degree 18 and is not homotopic to any injective virtual section because, by Theorem \ref{18}, there exists no multisection of degree 18. 
\end{rem}

Let us summarize this section: Multisections are in bijection with injective virtual sections (Corollary \ref{multi is inj virtual}). 
Homotopy as virtual sections and homotopy as multisections are equivalent conditions (Proposition \ref{mult vs virtual homotopy1} and \ref{mult vs virtual homotopy2}).

\section{Multisections and fundamental groups}
\label{fund group}
As we explained in the Introduction, our proof strategy is to relate the space $\x$ of smooth cubic plane curves and the classifying space $\mathrm{B}\slz$. We knew by the work of Dolgachev and Libgober \cite{DL} that $\pi_1(\x)$ fits into a short exact sequence:
$$0\to K\to\pi_1(\x)\to\slz\to 0$$
where $K$ is a finite Heisenberg group of order 27. Our goal in this section is to understand how $K$, which can be roughly viewed as the difference between $\x$ and $\mathrm{B}\slz$, interacts with multisections, and hereby to obtain some preliminary restrictions on multisections that will be used later in the proofs of main theorems.

\subsection{Fundamental group of $\x$: Dolgachev-Libgober's results}
In this subsection, we summarize some results of Dolgachev-Libgober \cite{DL} that we will use later.

Choose the Fermat cubic  $F(x,y,z)=x^3+y^3+z^3$ to be the base point of $\x$.  $C_F$ is homeomorphic to a torus. The torus bundle $\pi:\xp\to\x$ gives the following monodromy representation:
\begin{equation}
    \label{monodromy rep}
    \rho:\pix\to \mathrm{Aut}\ H_1(C_F;\ZZ)\cong \slz.
\end{equation}
Dolgachev-Libgober (Section 4 in \cite{DL}) proved that $\rho$ is surjective and gave the following explicit description of its kernel. Consider the subgroup of $\slc$ generated by the following two matrices:
\[
A=
\begin{bmatrix}
     0 & 0 & 1 \\
    1 & 0 & 0 \\
    0 & 1 & 0 \\    
\end{bmatrix}
,\ \ \  B=
\begin{bmatrix}
     1 & 0 & 1 \\
    0 & \zeta & 0 \\
    0 & 0 & \zeta^2 \\    
\end{bmatrix} 
, \ \text{where $\zeta:=e^{2\pi i/3}$}.\]
Both $A$ and $B$ have order 3 and their commutator is $\zeta I$ which is a generator of the center of $\slc$. Thus, the subgroup $K=\langle A,B\rangle$ is abstractly isomorphic to the Heisenberg group of order 27. 
\begin{rem}[\textbf{Notations for $\pgl$-actions}]
Since $\pgl$ acts naturally on $\CC P^2$, it also acts naturally on all of the following spaces: $\x, \xp, \xu, \xf.$ To avoid complicated notations, we will use the same notation $g\cdot x$ to denote the $\pgl$-actions on all of the spaces above. 
\end{rem}
Notice that $K$ preserves the Fermat cubic $F$. We have a well-defined orbit map 
\begin{align}
\label{orbit}    \mu:\slc/K &\to \x\\
\nonumber    g&\mapsto g\cdot F
\end{align}
inducing a map on the fundamental group (notice that  $\slc$ is 1-connected)
\begin{equation}
\label{K}
    \mu_*:K\cong\pi_1(\slc/K, I)\to \pix.
\end{equation}
Dolgachev-Libgober proved that $\mu_*$ is injective with image exactly the kernel of the monodromy representation $\rho$ in Section 4 of \cite{DL}. Thus, there is a short exact sequence:
\begin{equation}
    \label{DL}
    0\to K\xrightarrow{\mu_*} \pi_1(\x,F)\xrightarrow{\rho}  \slz\to0
\end{equation}
From this they computed the center of $\pi_1(\x)$.

\begin{prop}[Dolgachev-Libgober, Equation (4.9) in \cite{DL}]
\label{DL center}
$$Z(\pi_1(\x))=Z(K)\cong C_3.$$
\end{prop}

The action of the finite Heisenberg group $K$ on the Fermat cubic curve $C_F$ has the following simple description. The center $Z(K)$ acts trivially. The two generators $A$ and $B$ acts on the torus $C_F$ as translations of order 3 in two directions. Hence, 
\begin{equation}\label{translation}
    C_F/K \cong T^2.
\end{equation}
\subsection{Subgroups and quotients of $\pi_1(\x)$ and $\pi_1(\xp)$}
The goal of this subsection is to draw the commutative diagram in Proposition \ref{diagram proposition}. We will reach this goal by proving a sequence of lemmas. As in the previous subsection, we fix the Fermat cubic $F$ as the base point of $\x$.

\begin{lem}\label{pullback}
The orbit map $\slc/K\xrightarrow{\mu}\x$ in (\ref{orbit}) and the fiber bundle $\xp\xrightarrow{\pi}\x$ fit into the following commutative diagram
\begin{equation}\label{pullback bundle}
   \xymatrix{
\slc\times_K C_F \ar[r] \ar[d]_p & \xp \ar[d]^{\pi}\\
\slc/K \ar[r]^{\ \ \ \ \mu} & \x\\
} 
\end{equation}
where the map $\slc\times_K C_F\xrightarrow{p}\slc/K$ in the first column is given by the projection onto the first coordinate. Moreover, the fiber bundle $p$ is the pullback of the torus bundle $\pi$ via the map $\mu$.
\end{lem}
\begin{proof}
Recall that $\xp=\{(P,x)\in\x\times \CC P^2:x\in C_P\}$. 
By the definition of bundle pullback, the total space $\mu^*\xp$ of the pullback bundle is 
$$\mu^*\xp = \{(g,x)\in \slc/K\times \CC P^2 : x\in \mu(g)=g\cdot C_F\}.$$
To prove the lemma, it suffices to show that the following map gives an isomorphism of torus bundles over $\slc/K$
\begin{align}
\nonumber    \slc\times_K C_F&\longrightarrow \mu^*\xp\\
\label{pullback iso}    (g,y)&\longmapsto (g,g\cdot y).
\end{align}
This map is well-defined because for any $\gamma\in K$, the two points $(g,y)$ and $(g\gamma^{-1},\gamma y)$ are mapped to the same point $(g\gamma^{-1},g\gamma^{-1}\cdot \gamma y) = (g,g\cdot y)$ in $\mu^*\xp$. An inverse of the map (\ref{pullback iso}) is given by $(g,x)\mapsto (g,g^{-1} \cdot x)$.
\end{proof}

\begin{lem}
Applying $\pi_1$ to the diagram (\ref{pullback bundle}) with respect to appropriate choices of base points, we have the following commutative diagram of fundamental groups:
\begin{equation}
    \label{big diagram1}
    \xymatrix{
& & 0 \ar[d] & 0 \ar[d] & \\
0 \ar[r] & \ZZ^2 \ar[r] \ar[d] & \pi_1(\slc\times_K C_F) \ar[r] \ar[d] & K \ar[r] \ar[d]^{\mu_*} & 0\\
0 \ar[r] & \ZZ^2 \ar[r] & \pi_1(\xp) \ar[r]^{\pi_*} \ar[d] & \pi_1(\x) \ar[r] \ar[d]^\rho & 0\\
& & \slz \ar[d] \ar[r]^= & \slz \ar[d] & \\
& & 0 & 0 & \\
}
\end{equation}
\end{lem}
\begin{proof}
Notice that we suppress choices of base points in the diagram above. In fact, since all the spaces in the diagram (\ref{pullback bundle}) are path connected, different choices of base points will give us the same map up to conjugation. 

We have $\pi_1(\slc/K)\cong K$ because $\slc$ is simply connected. The top-right square commutes because the diagram (\ref{pullback bundle}) commutes. The second column on the right is exactly the short exact sequence (\ref{DL}) proven by  Dolgachev-Libgober. The first and the second rows in (\ref{big diagram1}) are obtained by applying the long exact sequence of homotopy groups to the two torus bundles in the commutative diagram in Lemma \ref{pullback}. The long exact sequences become short exact sequences because both base spaces have trivial $\pi_2$. That $\pi_2(\slc)=0$ is standard. We now show that $\pi_2(\x)=0$ using some classical results from algebraic geometry.

It is well-known that every smooth cubic plane curve is projectively equivalent to a curve in the \emph{Hesse canonical form}:
$$F_\lambda : x^3+y^3+z^3-3\lambda xyz=0$$
for some $\lambda\in \CC$. See \emph{e.g.} Lemma 1 in \cite{AD} for a proof. The family of curves $F_\lambda$ is called the \emph{Hesse pencil}. Moreover, $F_\lambda$ is smooth if and only if $\lambda^3\ne 1$. These classical results give us the following surjective map onto $\x$
\begin{align*}
    f: \pgl\times \Big(\CC\setminus\{1,\zeta,\zeta^2\}\Big)&\to \x\\
    (g,\lambda)&\mapsto g\cdot F_\lambda
\end{align*}
where $\zeta$ again denotes a primitive 3rd root of unity. Elements in $\pgl$ that preserve the Hesse pencil form a finite group of order 216, called the \emph{Hesse group}. The fibers of $f$ are precisely  orbits of the Hesse group acting diagonally and freely. Hence, $f$ is a covering map of order 216. We thus have
$$\pi_2(\x)\cong \pi_2(\pgl\times \CC\setminus\{1,\zeta,\zeta^2\})=0.$$
Again, we use the fact that $\pi_2(\pgl)\cong\pi_2(\slc)=0.$
\end{proof}

\begin{lem}\label{normal C3}
$\pi_1(\slc\times_K C_F)$ fits into the the following short exact sequence:
$$0\to Z(K)\to\pi_1(\slc\times_K C_F)\to \ZZ^2\to 0.$$
\end{lem}
\begin{proof}
Consider the map $\slc\times_K C_F\to C_F/K$ given by projection onto the second coordinate. We have  a fiber bundle
$$\xymatrix{
\slc/Z(K) \ar[r] & \slc\times_K C_F \ar[d]\\
&  \ \ \ \ \ \ C_F/K\cong T^2\\
}$$
Recall that $C_F/K\cong T^2$ by (\ref{translation}). Since $\pi_2(C_F/K)=0$, the fibration induces a short exact sequence
$$0\to \pi_1(\slc/Z(K))\to \pi_1(\slc\times_K C_F)\to \pi_1(T^2)\to 0.$$
Since $\slc$ is 1-connected, $\pi_1(\slc/Z(K))$ is naturally isomorphic to $Z(K)$. The lemma follows.
\end{proof}

\begin{lem}\label{normality}
  The fundamental groups of all the four spaces in (\ref{pullback bundle}) naturally contain $Z(K)$ as a normal subgroup.
\end{lem}
\begin{proof}
First of all, $\pi_1(\slc/K)\cong K$ contains $Z(K)$ as its center. $\pi_1(\x)$ also contains $Z(K)$ as its center, by Dolgachev and Libgober's result (Proposition \ref{DL center}). $\pi_1(\slc\times_KC_F)$ contains $Z(K)$ as a normal subgroup by Lemma \ref{normal C3}. All that is left is to check that $\pi_1(E)$ contains $Z(K)$ as a normal subgroup. This follows by a straightforward argument by chasing the commutative diagram (\ref{big diagram1}).
\end{proof}
\begin{rem}[\textbf{The geometric meaning of $Z(K)$ as an embedded subgroup}]\label{geometric meaning}
What does the normal subgroup $Z(K)$ mean or measure inside the fundamental groups of $\x$ and $\xp$ and other relevant spaces? Observe that the four spaces in the diagram (\ref{pullback bundle}) all have natural actions of $\pgl$ so that all of the maps in the diagram are $\pgl$-equivariant. Recall that $\pgl=\slc/Z(K)$. Hence, $Z(K)$ embeds inside each of the four $\pi_1$'s as the fundamental group of an $\pgl$-orbit with an appropriate choice of base points. Since $Z(K)$ is normal, different choices of base points give the same embedding. Hence, we can roughly view that $Z(K)$ measures those parts of the $\pi_1$'s of $\x$ and $\xp$ and other relevant spaces that come from the $\pgl$-action where $\pi_1(\pgl)=Z(K)$.
\end{rem}

Thanks to Lemma \ref{normality}, we can quotient out the common factor of $Z(K)$ from the relevant groups to simplify diagram (\ref{big diagram1}). Let's define
\begin{align}
\nonumber    &\ok:=K/Z(K)\cong (C_3)^2\\
\label{gamma}    &\gm:=\pi_1(\x)/Z(K)\\
\nonumber    &\gt:=\pi_1(\xp)/Z(K)
\end{align}
We also identify $\pi_1(\slc\times_K C_F)/Z(K)$ as $\ZZ^2$ thanks to Lemma \ref{normal C3}. Our discussion above gives the following proposition:
\begin{prop}\label{diagram proposition}
After taking those $Z(K)$-quotients as in (\ref{gamma}), the diagram (\ref{big diagram1}) becomes the following commutative diagram:
\begin{equation}
    \label{big diagram}
    \xymatrix{
& & 0 \ar[d] & 0 \ar[d] & \\
0 \ar[r] & \ZZ^2 \ar[r] \ar@{=}[d] & \ZZ^2 \ar[r] \ar[d] & (C_3)^2 \ar[r] \ar[d] & 0\\
0 \ar[r] & \ZZ^2 \ar[r] & \gt \ar[r]^{\pi_*} \ar[d] & \gm \ar[r] \ar[d]^\rho & 0\\
& & \slz \ar[d] \ar@{=}[r] & \slz \ar[d] & \\
& & 0 & 0 & \\
}
\end{equation}
\end{prop}
The diagram (\ref{big diagram}) can be considered  the diagram (\ref{big diagram1})  after ignoring those parts that come from the $\pgl$-actions, in the sense of Remark \ref{geometric meaning}. In the next subsection, we will show that the ignored factor $Z(K)$ is indeed irrelevant for our study of multisections (See Theorem \ref{C3 trivial} below).

\subsection{Fundamental group of multisection covers}
The goal of this subsection is to understand how the finite Heisenberg group $K$, which again can be viewed as the difference of $\x$ and $\mathrm{B}\slz$, interacts with multisections. Since $K$ fits into the following short exact sequence
$$0\longrightarrow Z(K)\longrightarrow K \longrightarrow \overline{K}\longrightarrow 0,$$
we will first understand how $Z(K)\cong C_3$ interacts with multisections and then understand how $\overline{K}\cong C_3\times C_3$ interacts with multisections. 

\begin{thm}
\label{C3 trivial}
If $\xn$ is a connected cover of $\x$ associated to a multisection, then $\pi_1(\xn)$ must contain $Z(K)$ as a normal subgroup. Here both $Z(K)$ and $\pi_1(\xn)$ are naturally identified as subgroups of $\pi_1(\x)$. 
\end{thm}
Notice that by Proposition \ref{DL center}, the subgroup $Z(K)$ is the center of $\pi_1(\x)$. Thus the statement $Z(K)\le \pi_1(\xn)$ makes sense without specifying a base point for $\xn$ even if $\pi_1(\xn)$ is not normal in $\pi_1(\x).$

\begin{rem}
One might (the authors once included) think that Theorem \ref{C3 trivial} is trivial by a simple ``proof" extending the idea in Remark \ref{geometric meaning}: $\pi_1(\xn)$ contains $Z(K)\cong \pi_1(\pgl)$ because the covering map $\xn\to\x$ is $\pgl$-equivariant. However, this argument fails because there is no natural $\pgl$-action on $\xn$ unless the multisection $\sigma:\x\to\xu$ is $\pgl$-equivariant, which is not true in general. The failure of a general multisection to be $\pgl$-equivariant is the main difficulty of Theorem \ref{C3 trivial}. Understanding such failure is also the main idea of the proof below. 
\end{rem}

\begin{proof}[Proof of Theorem \ref{C3 trivial}]
Again, let $F$ denote the Fermat cubic as in the previous subsections. Recall from (\ref{K}) that the map $K\hookrightarrow \pi_1(\x;F)$ is induced by the orbit map $\slc/K\to\x$. Thus, the embedding of $Z(K)$ in $\pix$ is induced by the orbit map 
\begin{align*}
    \slc/Z(K)=\pgl &\xrightarrow{\psi} \x\\
    g&\mapsto g\cdot F
\end{align*}
Hence, our claim that $\pi_1(\pgl)=Z(K)\le \pi_1(\xn)$ is equivalent to the following statement about map lifting:  there exists a map $\widetilde{\psi}$ making the  diagram commute:
\begin{equation}
    \label{lift a}
    \xymatrix{
& \xn \ar[d] \\
\pgl \ar[r]^\psi \ar@{-->}@/^5pt/[ru]^{\exists\widetilde{\psi}} & \x \\
}
\end{equation}
We will show that such a lift exists. First of all, let $\sigma:\x\to\xu$ denote the multisection that defines $\xn$. Let $n$ be its degree.  We consider the following map
\begin{align*}
    \pgl&\xrightarrow{\phi} \xu\\
    g&\mapsto g^{-1}\cdot \sigma(g\cdot F).
\end{align*}
If the multisection $\sigma$ is $\pgl$-equivariant, then $\phi$ would be a constant map. Hence, this map $\phi$ measures how far $\sigma$ deviates from being $\pgl$-equivariant based at the Fermat cubic $F$.  Moreover, $\sigma(g\cdot F)$ is a collection of $n$ points on the curve defined by $C_{g\cdot F}$, and thus $g^{-1}\cdot \sigma(g\cdot F)$ is a collection of $n$ points on $C_F$. In other words, the image of $\phi$ is entirely contained in a single fiber $\uconf_{n}(C_F)\subseteq \xu$. Let us now regard $\phi$ as a map $\pgl\to\uconf_{n}(C_F)$. Notice that $\uconf_n C_F$ is an Eilenberg-MacLane space (See Corollary 2.2 in \cite{FN}) and is also a finite dimensional manifold, and hence has a torsion-free  fundamental group. Thus, $\phi$ must induce a trivial map on fundamental groups. This implies that there exists a lift $\widetilde{\phi}$ making the following diagram commutes. 
\begin{equation}
    \label{lift b}
    \xymatrix{
& \uconf_{n-1,1} C_F \ar[d]^{n:1}\\
\pgl \ar[r]^\phi \ar@{-->}@/^10pt/[ru]^{\exists\widetilde{\phi}\ \ \ \ \ } & \uconf_n C_F \\
}
\end{equation}
Recall that $ \uconf_{n-1,1} C_F$ is the configuration space of $n$ unordered distinct points on $C_F$ with one of those $n$ points marked. Let $f:\uconf_{n-1,1} C_F\to C_F$ denote the ``forgetting" map that sends a marked configuration of $n$ points to its marked point. Now we can define the map $\widetilde{\psi}$ by 
\begin{align*}
   \widetilde{\psi}: \pgl&\longrightarrow\xn\\
    g&\longmapsto (g\cdot F,\ \  g\cdot f(\widetilde{\phi}(g)))
\end{align*}
Recall our assumption that $\xn$ is a cover associated to a multisection. Hence each point in $\xn$ can be uniquely represented as a pair $(P,y)$ such that $y$ is a point in $C_P$. 

We claim that $\widetilde{\psi}$ is a lift of $\psi$ making the diagram (\ref{lift a}) commute. First of all, we need to check that the image of $\widetilde{\psi}$ is in $\xn$. Indeed, since (\ref{lift b}) commutes, the point $f(\widetilde{\phi}(g))$ must belong to the configuration $\phi(g)= g^{-1}\cdot \sigma(g\cdot F)$. Hence, the point $g\cdot f(\widetilde{\phi}(g))$ must belong to the configuration $\sigma(g\cdot F)$. The  image of $\widetilde{\psi}$ is indeed in $\xn$. Finally, the diagram (\ref{lift a}) commutes because the first coordinate of $\widetilde{\psi}$ is exactly $\psi$.

So far we have proven that $Z(K)$ is a subgroup of $\pi_1(\xn)$. To show that $Z(K)$ is normal in $\pi_1(\xn)$, we just need to observe that $Z(K)$ is normal in the larger group $\pi_1(\x)$ by Proposition \ref{DL center}.
\end{proof}

Thanks to Theorem \ref{C3 trivial}, for any connected cover $\xn$ associated to a multisection of degree $n$, we will define 
\begin{align}
\label{gn}    &\gn:=\pi_1(\xn)/Z(K)
\end{align}
In particular, we define 
\begin{align}
\nonumber    &\gf:=\pi_1(\xf)/Z(K).
\end{align}
Recall that we defined $\ok:=K/Z(K)$ in (\ref{gamma}). 
\begin{thm}
\label{intersection is trivial}
If $\xn$ is a connected cover of $\x$ given by a multisection, then 
$$\ok\cap \gn=0.$$
\end{thm}
\begin{proof}
Consider the short exact sequence 
$$0\to\ZZ^2\to \gt\xrightarrow{\pi_*}\gm\to 0$$
in the middle row of (\ref{big diagram}).  Since $\gn$ comes from a multisection, its associated virtual section induces a group theoretic section $\gn\to\gt$ making the following diagram commute:
\begin{equation}
    \label{fin ind split}
    \xymatrix{
0 \ar[r] & \ZZ^2 \ar[r] & \gt \ar[r] & \gm \ar[r] & 0\\
& & & \gn \ar@{_{(}->}[u] \ar@{-->}[lu] & \\
}
\end{equation}
In other words, $\pi_*$ has a section when restricted to the finite index subgroup $\gn$ of $\gm$. Notice that the short exact sequence along the first row of (\ref{big diagram})
$$0\to\ZZ^2\to\ZZ^2\to \ok\to0$$
does not split even restricted to any nontrivial subgroup of $\ok\cong (C_3)^2$. Therefore, we must have that $\gn\cap \ok=0$. 
\end{proof}
Another way to formulate Theorem \ref{intersection is trivial} is that $K\cap \pi_1(\xn)=Z(K)$. 

Consider again the monodromy representation $\rho:\pi_1(\x)\to \slz$ of the universal cubic curve bundle $\pi$ as in (\ref{monodromy rep}). $\rho$ factors to a map (also called $\rho$ for simplicity)
$$\rho:\pi_1(\x)/Z(K)=\gm\to \slz$$ 
because $Z(K)=\langle\zeta I\rangle$ acts trivially on any cubic plane curve.
\begin{cor}
\label{9 mul}
Under the same assumption as in Theorem \ref{intersection is trivial}, the monodromy representation $\rho:\gm\to \slz$ maps the subgroup $\gn\le\gm$ isomorphically onto a  subgroup of $\slz$ of index $n/9$.
\end{cor}
\begin{proof}
Since $\ok$ is the kernel of the monodromy representation $\rho:\gm\to\slz$, the restriction $\rho|_{\gn}$ must be injective. To calculate the index, we have
$$n=[\pi_1\xn:\pi_1\x]=[\gn:\gm]=(\#\ker\rho)\cdot[\rho(\gn):\rho(\gm)]  = 9[\rho(\gn):\slz].$$
\end{proof}

As a consequence, any multisection must have degree $n$ such that $9| n$. This gives a new proof of Theorem \ref{old thm}.

We also obtain a new proof of the following result of Dolgachev-Libgober \cite{DL}.
\begin{cor}\label{flex iso slz}
The restriction of the monodromy representation $\rho:\gm\to\slz$ to the finite index subgroup $\gf\le\gm$ gives an isomorphism $\gf\cong \slz.$
\end{cor}
\begin{proof}
By Corollary \ref{9 mul}, the image of $\gf$ under $\rho$ is a subgroup of $\slz$ of index $9/9=1$ and hence is the entire $\slz$.
\end{proof}


\section{Multisections and cohomology}

In this section, we will prove one lemma about $H^1(\xn; \Z^2)$, where $\xn$ is the cover of $\x$ associated to a multisection and the $\ZZ^2$ coefficients are in the monodromy representation of the universal cubic curve bundle $\pi$ pulled-back to $\xn$.

We are interested in $H^1(\xn; \Z^2)$ because of this basic result in obstruction theory:

\begin{thm}[See \emph{e.g.} Theorem 6.12 of \cite{W}]
\label{representable}
Let $E \to B$ be a fiber bundle over a connected base. Suppose that the fiber is a $K(A,1)$ where $A$ is an abelian group. Then any two sections $s_1$ and $s_2$ are homotopic (relative to $B$) if and only if the corresponding obstruction class $\delta(s_1,s_2)$ is zero in $H^1(B;A)$. Here the coefficients of $H^1(B;A)$ are in the monodromy representation  $\pi_1(B)\to \mathrm{Aut}(A)$.

In particular, if $H^1(B;A)=0$, then sections of the bundle are unique up to homotopy or don't exist.
\end{thm}

We will apply Theorem \ref{representable} to prove Theorem \ref{9} and Theorem \ref{18} in the later sections. In order to do that, we need  to compute $H^1(\xn;\ZZ^2)$. Our next lemma reduces computing $H^1(\xn;\ZZ^2)$ to computing the group cohomology of $\gn$, which by Theorem \ref{intersection is trivial} is isomorphic to a finite index subgroups of $\slz.$

\begin{lem}
\label{H1 iso}
Let $\xn$ be a connected cover of $\x$ associated to a multisection. We have $H^1(\gn;\ZZ^2)\cong H^1(\xn;\ZZ^2)$.
\end{lem}
\begin{proof}
First of all, we have  $H^1(\xn;\ZZ^2)\cong H^1(\pi_1(\xn);\ZZ^2)$. It suffices to show that $H^1(\pi_1(\xn);\ZZ^2)\cong H^1(\gn;\ZZ^2)$. Since $\xn$ is connected, Theorem \ref{C3 trivial} tells us that there is a short exact sequence:
$$0\to Z(K) \to \pi_1(\xn)\to \gn\to0.$$
This short exact sequence induces a five-term exact sequence of cohomology groups with $\ZZ^2$ coefficients, whose initial terms are
$$0\to H^1(\gn; (\ZZ^2)^{Z(K)})\to H^1(\pi_1(\xn);\ZZ^2)\to H^1(Z(K);\ZZ^2)^{\gn}\to \cdots$$
The action of $Z(K)\cong C_3$ on $\ZZ^2$ is trivial because $Z(K)$ is generated by the matrix $\zeta I$, which acts trivially on any cubic plane curve. Hence, we have $(\ZZ^2)^{Z(K)}=\ZZ^2$ and $H^1(Z(K);\ZZ^2)=0$, which give an isomorphism $H^1(\pi_1(\xn);\ZZ^2)\cong H^1(\gn;\ZZ^2)$.
\end{proof}
Thanks to Theorem \ref{representable} and Lemma \ref{H1 iso}, in order to classify virtual sections, we need to understand $H^1(\gn;\ZZ^2)$ where $\gn$ is a subgroup of $\slz$ of index $n/9$ as in Corollary \ref{9 mul}.

\section{Proofs of Theorem \ref{9} and \ref{18}}

In this section, we show that multisections of low degrees are rigid -- we prove that multisections of degree 9 are unique up to homotopy (Theorem \ref{9}) and that multisections of degree 18 do not exist (Theorem \ref{18}).

\subsection{A general theorem}
First, we prove the following general theorem which we will use in the proofs of both Theorem \ref{9} and Theorme \ref{18}.
\begin{thm}
\label{H1 trivial}
Suppose a multisection induces a virtual section  $(\xn,p,s)$ such that $H^1(\xn;\ZZ^2)=0$. 
 Then we have the following results:
\begin{enumerate}
    \item Any other virtual section $s'$ from $\xn$ must be homotopic to $s$.
    \item  The cover $\xn$ factors through the flex cover $\xf$ defined in \emph{(\ref{def xf})}. In other words, there exists a map $f$ making the following diagram commute:
$$
\xymatrix{
\xn \ar@{-->}[rr]^{\exists\ f} \ar[rd]_p & & \xf \ar[ld]^{\pf}\\
& \x & \\
}$$
\item Let $\sflex:\xf\to\xp$ denote the flex virtual section defined by choosing the 9 flex points. Then $\sflex\circ f$ is homotopic to $s$ as virtual sections.
\end{enumerate}
\end{thm}

\begin{proof}
Applying Theorem \ref{representable} to the pullback of the torus bundle $\pi:\xp\to\x$ by the cover $p:\xn\to\x$,  we immediately obtain (1). 

If $\xn$ is disconnected, then part (2) will hold if we can prove the statement for each connected component of $\xn$. We can therefore assume $\xn$ to be connected without loss of generality. Since connected covers are classified by their fundamental groups,  in order to prove (2), it suffices to prove that $\pi_1(\xn)$ is conjugate into $\pi_1(\xf)$ as a subgroup of  $\pi_1(\x)$. Recall from Proposition \ref{C3 trivial} that $\pi_1(\x),\pi_1(\xn)$ and $\pi_1(\xf)$ all contain the normal subgroup $Z(K)\cong C_3$. Thus, it suffices to quotient out the common normal subgroup $Z(K)$ and to prove that $\gn$ is conjugate into $\gf$ as a subgroup of $\gm$. Recall the definitions of $\gn, \gf$ and $\gm$ in (\ref{gamma}) and (\ref{gn}). 

Consider the following short exact sequence, which is the first column of (\ref{big diagram}):
\begin{equation}
    \label{ses gamma}
\xymatrix{
0 \ar[r] & \ZZ^2 \ar[r] & \gt\ar[r] & \slz \ar[r] & 0.\\}
\end{equation}
The flex virtual section $\sflex: \xf\to\xp$ induces a map $\gf\to\gt$ which we will also denote by $\sflex$ for simplicity. The short exact sequence above has the following splitting:
$$
\xymatrix{
0 \ar[r] & \ZZ^2 \ar[r] & \gt\ar[r] & \slz \ar[r] & 0\\
& & & \gf \ar[u]_{\rho|_{\gf}}^\cong \ar[lu]^{\sflex} & \\
}
$$
$\rho|_{\gf}$ is an isomorphism by Corollary \ref{flex iso slz}. The composition $\sflex\circ(\rho|_{\gf})^{-1}$ gives a group-theoretic section of the short exact sequence. In a similar way, the virtual section $s$ gives another group-theoretic section of the short exact sequence (\ref{ses gamma}) defined only on a finite index subgroup of $\slz$. More precisely, we have:
$$
\xymatrix{
0 \ar[r] & \ZZ^2 \ar[r] & \gt\ar[r] & \slz \ar[r] & 0\\
& & & \gn \ar@{_{(}->}[u]_{\rho|_{\gn}} \ar[lu]^{s} & \\
}
$$
The injectivity of  $\rho|_{\gn}$ follows from Corollary \ref{9 mul}. In summary, the two sections constructed above give two sections of the following extensions of $\rho(\gn)\le \slz$:
$$
\xymatrix{
0 \ar[r] & \ZZ^2 \ar[r] & \ZZ^2\rtimes \rho(\gn) \ar[r] & \rho(\gn)\ar@{-->}@<-0.5ex>@/^-10pt/[l]_{\ \ \ \  \ \ s\ \ \ \ } \ar@{-->}@<0.5ex>@/^10pt/[l]^{\sflex} \ar[r] & 0.}
$$
By Lemma \ref{H1 iso}, we have  $H^1(\gn;\ZZ^2)\cong H^1(\xn;\ZZ^2)=0$. This implies that sections of the short exact sequence above are unique up to conjugation. In particular, $s$ and $ \sflex$ are conjugate: There exists an element $\tilde{\gamma}\in\gt$ such that
\begin{equation}
    \label{conjugate}
    \tilde{\gamma}\ s\Big(\rho(\gn)\Big)\tilde{\gamma}^{-1} =\sflex \Big(\rho(\gn)\Big).
\end{equation}
We now apply the map $\pi_*$ which is one of the horizontal maps in (\ref{big diagram}) to both sides of (\ref{conjugate}). Let $\gamma:=\pi_*(\tilde{\gamma})\in\gm$. We have
$$    {\gamma}\ \pi_*s\Big(\rho(\gn)\Big){\gamma}^{-1} =\pi_*\sflex \Big(\rho(\gn)\Big).
$$
On the left hand side, we have
$$LHS= \gamma \ \gn \gamma^{-1}.$$
On the right hand side, we have
$$RHS \subseteq \pi_*\sflex(\slz) = \gf$$
Therefore, we have $\gamma \ \gn \gamma^{-1}\subseteq \gf$ and part (2) follows.

Part (3) follows immediately from (1) and (2). 
\end{proof}

\subsection{Proof of Theorem \ref{9}: Multisections of degree 9 are unique}
The proof of Theorem \ref{9} will  be short now that we have proven the general Theorem \ref{H1 trivial}. Suppose that $\sigma$ is any multisection of degree 9. Our goal is to show that $\sigma$ is homotopic to the multisection defined by picking the 9 flex points. By Proposition \ref{mult vs virtual homotopy2}, it suffices to to show that $\sigma$ induces a virtual section that is homotopic (as virtual sections) to the flex virtual section. Indeed, we will have:
\begin{prop}
\label{no 9}
Any injective virtual section $(\xn,p,s)$ of degree $9$ is homotopic via virtual sections to the flex virtual section.
\end{prop}
\begin{proof}
We first claim that $\xn$ must be connected. Let $Y$ be a connected component of $\xn$. By Corollary \ref{9 mul}, we have $\deg(Y/X)\ge 9$. However, we also have $\deg(Y/X)\le \deg(\xn/X)= 9$. So $Y$ must be the entire $\xn$.

The claim will follow from Theorem \ref{H1 trivial} if we can show that $H^1(\xn;\ZZ^2)=0$, which by Lemma \ref{H1 iso} is equivalent to showing that $H^1(\gm_9;\ZZ^2)=0$ since $\xn$ is connected. We know that $\gm_9$ is isomorphic to $\slz$ by Corollary \ref{9 mul}. It is well known that $H^1(\slz;\ZZ^2)=0$ (see Theorem 3.7, part (3) in \cite{CCS}). The proposition now follows from Theorem \ref{H1 trivial}.
\end{proof}

\subsection{Proof of Theorem \ref{18}: Multisections of degree 18 do not exist}
Suppose that $\sigma$ is a multisection of degree 18. Let $(\xn,p,s)$ denote the induced virtual section. The proof will proceed in three steps. 

\noindent\textbf{Step 1: Show that $\pmb{H^1(\xn;\Z^2) =0}$.}  We first prove that $H^1(\xn;\Z^2) =0$ in order to apply Theorem \ref{H1 trivial}.

It is possible that $\xn$ is a disconnected cover. By Corollary \ref{9 mul}, each connected component of $\xn$ is itself a cover of $\x$ with degree a multiple of 9. There are  two possibilities: either (1) $\xn$ is connected or (2) $\xn$ has two connected component, each of which is a degree 9 cover of $\x$. We claim that in either case we have $H^1(\xn;\ZZ^2)=0.$ 

If (2) is the case, then by the proof of Theorem \ref{no 9} above, each connected component $Y$ of $\xn$ has $H^1(Y;\ZZ^2)=0$. Thus, $H^1(\xn;\ZZ^2)=0$.

Let's consider case (1). By Lemma \ref{H1 iso}, if $\xn$ is connected, then $H^1(\xn;\Z^2) = H^1 (\gm_{18};\Z^2) $. Here $\gm_{18}$ is the unique index-2 subgroup of $SL_2(\Z)$ by Corollary \ref{9 mul}.
\begin{lem}
$H^1(\gm_{18};\ZZ^2)=0.$
\end{lem}
\begin{proof}
First, we prove that all elements in $H^1(\gm_{18};\ZZ^2)$ are 2-torsion. Recall that 
$$H^1(\gm_{18};\ZZ^2) = \frac{Z(\gm_{18};\ZZ^2)}{B(\gm_{18};\ZZ^2)}$$
where $Z(\gm_{18};\ZZ^2)$ is the group of cocylces, which are functions $\gm_{18}\to \ZZ^2$ satisfying the cocycle conditions, and $B(\gm_{18};\ZZ^2)$ is the group of coboundaries. Since $\gm_{18}$ is the (unique) index-2 subgroup of $\slz$, we have $-I\in \gm_{18}$. Suppose $\phi$ is any cocyle in $Z(\gm_{18};\ZZ^2)$, if $\phi(-I)=0$, then for any $g\in Z(\gm_{18};\ZZ^2)$, the cocyle conditions give:
$$\phi(-g)= \phi((-I)g)=\phi(-I)-\phi(g)=-\phi(g)\in \ZZ^2$$
$$\phi(-g)= \phi(g(-I))=\phi(g)+g\phi(-I)=\phi(g)\in \ZZ^2$$
Hence, $\phi$ must be the zero cocycle. More formally, this implies that the map $Z^2(\gm_{18}, \ZZ^2) \to \ZZ^2 $ defined by  $\phi \mapsto \phi(-I)$ is injective. Moreover, we claim that under this embedding, we have $2\ZZ^2 \subseteq B^2(\gm_{18}, \ZZ^2) \subseteq Z^2(\gm_{18}, \ZZ^2) \subseteq \ZZ^2$. Indeed, for any element $2v\in 2\ZZ^2$, the cocycle $\phi_v$ defined by $\phi_v(g)= gv-v$ is a coboundary such that $\phi_v(-I)=-2v$. Hence, $2\ZZ^2\subseteq B^2(\gm_{18}, \ZZ^2)$. Therefore, all elements in $H^1(\gm_{18};\ZZ^2)$ must be 2-torsion. 

Finally, we show that $H^1(\gm_{18};\ZZ^2)$ must be zero. Suppose not. Take any cocycle $\phi$ that is not a coboundary. By the previous paragraph, we know that $2\phi$ must be a coboundary. So there exists an element $v\in\ZZ^2$ such that $2\phi(g)=gv-v$ for any $g\in \gm_{18}.$ Moreover, since $\phi$ is not a coboundary, the element $v$ is not in $2\ZZ^2.$ Now we have
$$\forall g\in \gm_{18}, \ \ \ \ gv-v=2\phi(g)=0 \text{ (mod 2)}.$$
Let $\overline{v}$ denote that image of $v$ in $(\ZZ/2\ZZ)^2$. Consider the action of $\slz$ on $(\ZZ/2\ZZ)^2$. Let $\gm_{\overline{v}}$ denote the stabilizer of the nonzero element $\overline{v}$. The statement above says that $\gm_{18}\le \gm_{\overline{v}}$. The index of the subgroup $\gm_{\overline{v}}$ in $\slz$ is 3 because $\slz$ acts transitively on the set of three nonzero elements in $(\ZZ/2\ZZ)^2$. However, $\gm_{18}$ is an index-2 subgroup of $\slz$. We reach a contradiction.
\end{proof}


\noindent\textbf{Step 2: Construct a vector field.}
In this step, we will construct a nonvanishing vector field of a vector bundle. In the next and final step, we will prove that such a vector field cannot exist because of a nonzero obstruction classes, finishing the proof  by contradiction.


In the previous step, we proved that $H^1(\xn;\ZZ^2)=0$. By Theorem \ref{H1 trivial}, there is a covering map $f:\xn\to\xf$ such that the virtual section $s$ is homotopic to $\sflex\circ f.$ Let $s_t: \xn\to \xp$ be such a homotopy  that $s_0=\sflex\circ f$ and $s_1=s.$ 
$$
\xymatrix{
\xn\  \ar@{^{(}->}[r]^s \ar[d]^f_{2:1} & \xp\\
\xf \ar[ru]_{\sflex} & \\
}
$$

\begin{figure}[ht]
    \centering
    \includegraphics[scale=.5]{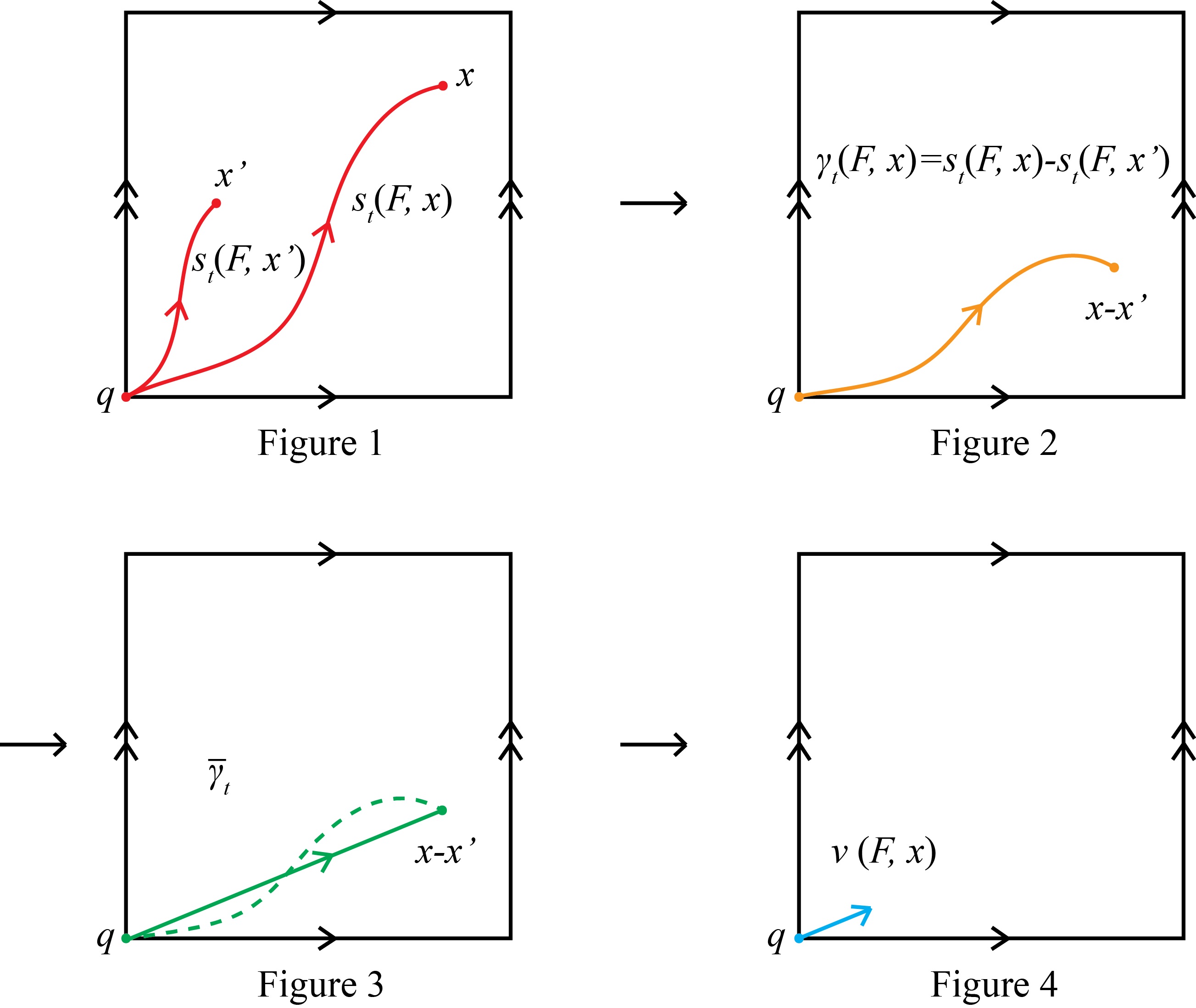}
     \caption{An illustration of the 4-step construction of the vector $v(F,x)$ out of the homotopy $s_t(F,x)$ on a fixed cubic curve $C_F$.}
     \label{construct vector}
\end{figure}

Let $\xi$ denote the (real) rank-2 vector bundle over $\xf$ whose fiber over a point $(F,q)$ is the  tangent space $T_qC_F$ of the curve $C_F$ at its flex point $q$. Let $f^*\xi$ denote the pullback of $\xi$ over $\xn$ by the map $f:\xn\to\xf.$ Our goal in this step is to construct a nonvanishing vector field of the bundle $f^*\xi$ using the homotopy $s_t$. We will describe the construction in the next paragraph. See Figure \ref{construct vector} for pictures of the constructions.

Since the cover $\xn$ of $\x$ comes from a multisection, each point on $\xn$ can be uniquely represented as a pair $(F,x)$ where $x$ is a point on $C_F$. Moreover, the 2:1 cover $f:\xn\to\xf$ maps two distinct points $(F,x), (F,x')\in \xn$ onto the same $(F,q)\in \xf$. For now, we will consider a fixed curve $C_F$ and a given point $x$ in it. For simplicity, we will just write $f(x)=q$ instead of $f(F,x)=(F,q)$. There is a unique point $(F,x')\in\xn$ such that $f(x)=f(x')=q$ but $x\ne x'$. Define a path on $C_F$ by
\begin{align*}
\gamma_t{(F,x)}
    :=s_t(F,x)-s_t(F,x') && t\in[0,1]
\end{align*} 
where the group law is taken in the elliptic curve $C_F$ with $q=f(x)$ as the identity. $\gamma_t{(F,x)}$ is a continuous path on $F$ which starts at 
\begin{align*}
    s_0(F,x)-s_0(F,x')&=\sflex\circ f(F,x)- \sflex\circ f(F,x')\\
    &= \sflex (F,f(x)) -\sflex (F,f(x'))&\text{note $f(x)=f(x')=q$}\\
    &= q-q=q &\text{since $q$ is the identity}
\end{align*}
and ends at 
$$s_1(F,x)-s_1(F,x') = s(F,x)-s(F,x')=x-x'\ne q.$$
The last equality of the equation above follows from our assumption that $x\ne x'$ so their difference is not equal to the identity $q$ under the group law. The path $\gamma_t{(F,x)}$ has two distinct endpoints. Equip $C_F$ with its unique flat Riemannian metric compatible with its complex structure with unit area. Every path on the flat torus $C_F$ with distinct endpoints is homotopic relative to its endpoints to a unique unit-speed geodesic by a linear homotopy. Let $\overline{\gamma}_t{(F,x)}$ be the unique unit-speed geodesic homotopic to $\gamma_t{(F,x)}$. Consider its tangent vector at the start point
$$v{(F,x)}:=\frac{d}{dt}\overline{\gamma}_t{(F,x)}(t)\Bigr|_{t=0}\in T_{q}C_F.$$
In particular, $v{(F,x)}\ne0$. 

The vector $v{(F,x)}$ varies continuously with $(F,x)\in\xn$ because of the following facts: $s_t$ is a homotopy of virtual sections. $f$ is also a continuous map of the parameter spaces $\xn\to\xf$. The elliptic curve and the Riemannian metric on $C_F$ vary continuously with $F\in\x$. 
To sum up, the map 
\begin{align*}
\xn&\longrightarrow T_qC_F\\
(F,x)&\longmapsto v{(F,x)}
\end{align*}
is a continuous nonvanishing vector field of the bundle $f^*\xi$.

\noindent\textbf{Step 3: Computing the Euler class of $\xi$.} Since the rank-2 vector bundle $f^*\xi$ has a nonvanishing section, its Euler class 
$$e(f^*\xi)\in H^2(\xn;\ZZ)$$
must be zero. Now consider the composition of the following two maps
\begin{align*}
    H^2(\xf;\ZZ)&\xrightarrow{f^*} H^2(\xn;\ZZ)\xrightarrow{f^!} H^2(\xf;\ZZ)
\end{align*}
where $f^!$ is the transfer homomorphism associated to the 2-to-1 cover $f:\xn\to\xf$. The composition $f^!\circ f^*$ is multiplication by $\deg f=2$. Hence, we have
$$2e(\xi)=f^!\circ f^* \Big(e(\xi)\Big) = f^! e(f^*\xi) =f^!(0)=0.$$
Therefore, $e(\xi)$ must be 2-torsion in $H^2(\xf;\ZZ)$. We will reach a contradiction once we  prove the following proposition. 
\begin{prop}\label{torsion 12}
$e(\xi)$ is a torsion element of order 12 in $H^2(\xf;\ZZ)$.
\end{prop}
\begin{proof}
We will prove this proposition by proving a sequence of lemmas.

Recall that every orientable rank-2 vector bundle is a pullback of the universal bundle over $B\slr$ by a classifying map unique up to homotopy. Let $\phi_\xi:\xf\to B\slr$ denote the classifying map of $\xi$. In the following lemma, we claim that $\phi_\xi$ is the composition of two  maps that we understand relatively well. 

Let $\eta$ denote the torus bundle over $\xf$ whose fiber over $(F,q)\in\xf$ is $C_F$. Equivalently, $\eta$ is the pullback of the universal cubic curve $\pi:\xp\to\x$ by the flex cover $\xf\to\x$.
\begin{lem}\label{BSL}
The classifying map $\phi_\xi:\xf\to B\slr$ of the bundle $\xi$ is the composition of the following two maps up to homotopy
    $$\xf\xrightarrow[]{\psi} \mathrm{B}\slz\xrightarrow[]{\iota} B\slr$$
    where $\iota$ is induced by the natural inclusion $\slz\hookrightarrow\slr$ and $\psi$ induces the monodromy representation $\pi_1(\xf)\to\slz$ of the torus bundle $\eta$ over $\xf$.
\end{lem}
\begin{proof}
First, we claim that $\phi_\xi$ is the composition of two maps in the commutative diagram (up to homotopy) below:
\begin{equation}
    \label{classifying maps}
    \xymatrix{
\xf \ar[rr]^{\phi_\xi} \ar[rd]_{\phi_\eta} & & B\slr\\
& B\mathrm{Diffeo}^+(T^2,q) \ar[ru]_{\tau} & \\
}
\end{equation}
In the diagram above, $\mathrm{Diffeo}^+(T^2,q)$ denotes the group of orientation-preserving diffeomorphisms of the torus $T^2$ fixing a base point $q\in T^2$. The map $\phi_\eta$ denotes the classifying map of the smooth orientable torus bundle $\eta$. The map $\tau$ is induced by the map $\mathrm{Diffeo}^+(T^2,q)\to \slr$ obtained by taking the derivative of a diffeomorphism at the fixed point $q$. The diagram commutes up to homotopy because $\xi$ is the vector bundle whose fiber over a point $(F,q)$ is the tangent vector space $T_q(C_F)$.

Thanks to a result of Earle-Eells (Theorem 1 in \cite{EE}), we know that  $\mathrm{Diffeo}^+_0(T^2,q)$, the component of the identity in $\mathrm{Diffeo}^+(T^2,q)$, is contractible. The quotient $\mathrm{Diffeo}^+(T^2,q)/\mathrm{Diffeo}^+_0(T^2,q)$ is the mapping class group of the torus, which is isomorphic to $\slz$ by taking its action on $H_1(T^2;\ZZ)\cong \ZZ^2$ (See \emph{e.g.} Theorem 2.5 in \cite{FM} for a proof of this well-known fact). Hence,  we have a homotopy equivalence
$$B\mathrm{Diffeo}^+(T^2,q) \xrightarrow{\sim} \mathrm{B}\slz.$$
We can thus replace (\ref{classifying maps}) by the following commutative diagram (up to homotopy):
\begin{equation}
    \label{classifying maps2}
    \xymatrix{
\xf \ar[rr]^{\phi_\xi} \ar[rd]_{\psi} & & B\slr\\
& \mathrm{B}\slz \ar[ru]_{\iota} & \\
}
\end{equation}
By our construction above, $\psi_*:\pi_1(\xf)\to\slz$ is exactly the monodromy representation of the torus bundle $\eta$. Moreover, the composition
$$\slz\to \mathrm{Diffeo}^+(T^2,q)\to \slr$$
is the standard inclusion up to conjugation. Hence, $\iota$ is induced by the natural inclusion $\slz\hookrightarrow\slr$. 
\end{proof}

Therefore, in order to understand what $\phi_\xi$ does on $H^2$, it suffices to understand what $\psi$ and $\iota$ do on $H^2$, respectively.

\begin{lem}
\label{pi inj}
$\psi$ induces an injective map $\psi^*: H^2(\mathrm{B}\slz;\ZZ)\hookrightarrow H^2(\xf;\ZZ)$.
\end{lem}

\begin{proof}

Recall that by Corollary \ref{flex iso slz}, we have  
$$\gf:=\pi_1(\xf)/Z(K)\cong \slz$$
where the isomorphism is given by the restriction of the monodromy representation $\rho:\pi_1(\x)\to \slz$ to the subgroup $\pi_1(\xf)$, or equivalently, given by the monodromy representation of the pullback bundle $\eta$.

Let $Y$ be the cover of $\xf$ associated to the subgroup $Z(K)$ in $\pi_1(\xf)$. So we have that $\pi_1(Y)=Z(K)\cong C_3$ and that the deck group of the cover $Y\to\xf$ is $\slz$. This cover induces the following  bundle over $\mathrm{B}\slz$ with fiber $Y$:
$$
\xymatrix{
Y \ar[r] & Y\times_{\slz} E\slz \ar[d] \\
& \mathrm{B}\slz & \\
}
$$
Since $E\slz$ is contractible, the total space of the bundle is homotopy equivalent to $Y/\slz = \xf$. In other words, we have a homotopy fibration
\begin{equation}
    \label{fibration}
    Y\to \xf\xrightarrow{\psi} \mathrm{B}\slz.
\end{equation}
The projection map is homotopic to $\psi$ because, by Lemma \ref{BSL}, the two maps induced the same map on fundamental groups up to conjugation. Let us just identify the projection map as $\psi$ since homotopy makes no difference in the discussion below.

We will prove the lemma by analyzing the Serre spectral sequence of the fibration (\ref{fibration}):
$$E_2^{p,q}= H^p(\mathrm{B}\slz;H^q(Y;\ZZ))\ \  \Longrightarrow\ \  H^{p+q}(\xf;\ZZ)$$
In particular, we have the following commutative diagram of the edge morphism:
\begin{equation}
    \label{edge homo}
    \xymatrix{
H^2(\mathrm{B}\slz;\ZZ) \ar[r]^{\psi^*} \ar@{=}[d] & H^2(\xf;\ZZ)\\
E^{2,0}_2 \ar@{->>}[r] & E^{2,0}_\infty \ar@{^{(}->}[u]\\
}
\end{equation}
On the $E_2$-page, the only possibly nontrivial differential into $E^{2,0}_2$ is from $E^{0,1}_2 = H^0(\mathrm{B}\slz;H^1(Y;\ZZ))$. However, $H^1(Y;\ZZ) = \mathrm{Hom}(\pi_1Y,\ZZ) =\mathrm{Hom}(C_3,\ZZ)=0$. Hence, the surjection $E^{2,0}_2\twoheadrightarrow E^{2,0}_\infty$ is an isomorphism. From the diagram (\ref{edge homo}), we can see that $\psi^*:H^2(\mathrm{B}\slz;\ZZ)\to H^2(\xf;\ZZ)$ is injective as we claimed. 
\end{proof}

\begin{lem}\label{iota}
$\iota$ induces a surjection $\iota^*: H^2(B\slr;\ZZ)\to H^2(\mathrm{B}\slz;\ZZ)$.
\end{lem}

\begin{proof}
We know that 
$$B\slr\simeq B\mathrm{U}(1)\simeq \CC P^\infty.$$
We will prove the lemma by restricting $\iota^*$ to two finite cyclic subgroups of $\slz$. Define
$$A:=\begin{bmatrix}
0 & -1 \\
1 & 0 
\end{bmatrix},\ \ \ 
B:=\begin{bmatrix}
0 & -1 \\
1 & 1 
\end{bmatrix}.
$$
$A$ has order 4 and $B$ has order 6. They together generate $\slz$. Consider the map
$$\langle A\rangle\cong C_4\hookrightarrow \slz\hookrightarrow \slr.$$
Observe $A$ is an element in $SO(2) = \mathrm{U}(1)$ which acts on the circle by $\pi/2$-rotation. A model for $B\mathrm{U}(1)$ is
$$B\mathrm{U}(1) = (\CC^\infty\setminus \{0\})/\mathrm{U}(1) \simeq \CC P^\infty$$
A model for $BC_4$ is
$$BC_4 = (\CC^\infty\setminus \{0\})/\langle A\rangle \simeq S^\infty/\langle A\rangle = L(4)$$
where $L(4)$ is an infinite dimensional lens space. Using these explicit models, it is straightforward to check that the inclusion $\langle A\rangle$ induces a surjective map on $H^2(B \mathrm{U}(1);\ZZ)\cong \ZZ\to H^2(B C_4;\ZZ)\cong \ZZ/4\ZZ$.

If we replace $A$ with $B$, the exact same argument does not apply verbatim because $B$ is no longer an element in $SO(2)$. However, $B$ is conjugate in $\slr$ to an  element $B'$ in $SO(2)$ which acts as the $\pi/3$-rotation of the circle. Therefore, the inclusion $C_6=\langle B\rangle \hookrightarrow \slr$ is homotopic to the inclusion $C_6=\langle B'\rangle \hookrightarrow SO(2)\subseteq \slr$. Applying the  argument above to $B'$, we conclude that the inclusion $\langle B\rangle$ induces a surjective map on $H^2(B \mathrm{U}(1);\ZZ)\cong \ZZ\to H^2(B C_6;\ZZ)\cong \ZZ/6\ZZ$.

To sum up, the map $\iota^*$ is surjective when restricted to $H^2(BC_4 ;\ZZ)\cong \ZZ/4\ZZ$ and  $H^2(BC_6 ;\ZZ)\cong \ZZ/6\ZZ$. Therefore, $\iota^*:H^2(B\slr;\ZZ)\cong \ZZ\to H^2(\mathrm{B}\slz;\ZZ)\cong \ZZ/12\ZZ$ must be surjective. 
\end{proof}

Finally, we now finish the proof of Proposition \ref{torsion 12}. The \emph{universal Euler class} is a generator $e$ of $H^2(B\slr;\ZZ)\cong \ZZ$. The Euler class of the rank-2 vector bundle $\xi$ is the pullback of $e$ by its classifying map $\phi_\xi$. Apply the $H^2$ functor to diagram (\ref{classifying maps2}), we have
$$e(\xi)=\phi_\xi^* (e) = \pi^*\iota^*(e).$$
Since $\iota^*$ is surjective by Lemma \ref{iota}, $\iota^*(e)$ is a generator of $H^2(\mathrm{B}\slz;\ZZ)\cong \ZZ/12\ZZ$. Since  $\pi^*$ is injective by Lemma \ref{pi inj}, $\pi^*\iota^*(e)$ remains an element of order 12 in  $H^2(\xf;\ZZ).$ Our proof of Proposition \ref{torsion 12} is complete.
\end{proof}
We now reach a contradiction in the three steps above. Hence, the multisection $\sigma$ of degree 18 cannot exist. \qed

\section{Algebraic constructions of multisections}\label{alg cons}
In this section, we describe how classical constructions from algebraic geometry give us multisections of various degrees as in Theorem \ref{old AG}. More importantly, we compute the monodromy group associated to those algebraic multisections, which will allow us to distinguish these algebraic constructions from the topological constructions which we will introduce in the next section. 

\subsection{Algebraic constructions}

Every smooth cubic plane curve $C$ has 
\begin{itemize}
\item 9 flex points where the tangent line intersect $C$ with  multiplicity $3$ (first attributed to Maclaurin; see Introduction of \cite{AD} for a brief history), and
\item 27 sextatic points where the osculating conic intersects $C$ with multiplicity $6$ (Cayley \cite{Cayley}).
\end{itemize}
More generally, Gattazo \cite{g} made the following definition:
\begin{defn}[Gattazzo]
\label{type 3m}
A {point of type $3m$} on a smooth cubic plane curve $C$ is where $C$ intersects an irreducible cure of degree $m$ with multiplicity $3m$. Equivalently, if we choose an arbitrary flex point as the identity for $C$, then a {point of type $3m$} is precisely a $3m$-torsion point that is not a $3k$-torsion point for any $k<m$.
\end{defn}
The equivalence of the two definitions above follows from this fact: $3m$ points $P_i$ for $i=1,...,3m$ on a cubic curve $C$ are on another curve of order $m$ if and only if $\sum_{i-1}^{3m} P_i=0$ on the elliptic curve $C$ with a flex point chosen as identity.  See page 392 in \cite{c} for a proof of this fact. Under this definition, a flex point is of type 3, while a sextatic point is  of type 6. 

The \emph{Jordan 2-totient function} $J_2$ can be defined as: 
\begin{equation}
\label{number of points of type 3m}
    J_2(m) = m^2\prod_{p|m,\ p \text{ prime}}(1-\frac{1}{p^2}).
\end{equation}
$J_2(m)$ is equal to the number of torsion points on an elliptic curve with order $m$. The number of points of type $3m$ on a smooth cubic plane curve is $9J_2(m)$.

\subsection{Multisections from torsion constructions}
The map 
\begin{align}
\label{sigma_m}    \sigma_m:\x&\longrightarrow \xu\\
\nonumber    F&\longmapsto \{\text{points of type $3m$ on $C_F$}\}
\end{align}
defines a multisection of degree $9J_2(m)$. Let $\xm$ denote the associated cover of $\x$:
\begin{equation}
\label{xm define}
    \xm:=\{(F,x)\ |\ x \text{ is a point of type $3m$ on $C_F$}\}. 
\end{equation}

Recall that $\rho:\pi_1(\x)\to\slz$ denotes the monodromy representation  of the universal cubic curve $\pi:\xp\to\x$ as in (\ref{monodromy rep}). 

\begin{lem}
\label{conjugate to congruence subgroup}
For each positive integer $m$, the cover $\xm$ is connected. Moreover, under the monodromy representation $\rho:\pi_1(\x)\to\slz$, the image of the subgroup $\pi_1(\xm)$  is conjugate to the following congruence subgroup in $\slz$ 
\begin{equation}
    \label{cong subg}
\Gamma_1(m)=\bigg\{g\in\slz \ |\ g=
\begin{bmatrix}
1&*\\
0&1
\end{bmatrix}
\mod m\bigg\}.
\end{equation}
\end{lem}
\begin{proof}

Let $Y$ be any connected component of $\xm$. We will show that the image of $\pi_1(Y)$ under $\rho$ is conjugate to a subgroup of $\Gamma_1(m)$ in $\slz$. Next, by comparing the indices, we show that $Y$ must be the entire $\xm$ and that $\rho(\pi_1(\xm))$ must be conjugate to the entire $\Gamma_1(m)$.

Take any $(F,x)$ in $Y$. Let $z$ denote the element $3q$ where $q$ is an arbitrary flex point on $C_F$ and the group law is with $x$ as the identity. We first claim that $z$ is well-defined. Indeed, if $p$ is any other flex point, then $p-q$ is a 3-torsion on $(C_F,p)$ and hence is also a 3-torsion on $(C_F,x)$ because the difference between two points remains unchanged when we choose a different identity of an elliptic curve. Moreover, we have $mz=3mq=0$ on $(C_F,x)$ because $q-x$ is a $3m$-torsion on $(C_F,q)$ and hence also a $3m$-torsion on $(C_F,x)$. Therefore, $z$ is an $m$-torsion point on $(C_F,x)$. Furthermore, $z$ is of order exactly $m$ since $x$ is a point of type $3m$.

By the construction above, every element $(F,x)\in Y$ comes equipped with a marked element $z=z(F,x)$ on $C_F$ which is a torsion point of order $m$ on the elliptic curve $(C_F,x)$. Consequently, the monodromy action of $\pi_1(Y)$ on $H_1(C_F;\ZZ/m\ZZ)$ must fix an element $v$ that corresponds to the $m$-torsion point $z$. Moreover, since $z$ is of order $m$, by a change of basis we can choose $v$ to be $e_1=[1,0]^\mathrm{T}$ in $\mathbb{F}_3^2$. Hence, $\rho(\pi_1(Y))$ is conjugate to a subgroup of $\Gamma_1(m)$ in $\slz$.

Finally, the degree of the covering $Y$ over $\x$ is
\begin{align*}
&\deg(Y/\x) = [\pi_1(Y):\pi_1(\x)]\\
&= [\pi_1(Y)/Z(K):\pi_1(\x)/Z(K)]&\text{by Theorem \ref{C3 trivial}}\\
&= 9\cdot [\rho(\pi_1(Y)): \slz]&\text{by Corollary \ref{9 mul}}\\
&\ge 9\cdot [\Gamma_1(m): \slz]]&\text{$\rho(\pi_1(Y))$ is conjugate to a subgroup of $\Gamma_1(m)$}\\
&= 9\cdot J_2(m)
\end{align*}
On the other hand, $\deg(Y/\xn)\le \deg(\xm/\x)=9\cdot J_2(m)$ because $Y$ is a component of the cover $\xm$. Hence, both inequalities above must be equalities. We have that $\rho(\pi_1(Y))$ is conjugate to the entire  $\Gamma_1(m)$ and that $\xm=Y$ is connected.
\end{proof}

\begin{cor}
\label{cong is free}
When $m\ge 4$, the group 
$$\gn=\pi_1(\xm)/Z(K)$$
as defined in (\ref{gn}) is free. Here the degree $n$ is $9J_2(m)$.
\end{cor}
\begin{proof}
By Theorem \ref{intersection is trivial}, we know that $\rho$ maps $\gn$ isomorphically onto a subgroup of $\slz$ which, by Lemma \ref{conjugate to congruence subgroup}, is conjugate to $\Gamma_1(m)$ as in (\ref{cong subg}). It now takes a standard argument to check that $\Gamma_1(m)$ is free when $m\ge4$. It suffices to show that  $\Gamma_1(m)$ is torsion-free when $m\ge4$ because $\slz$ acts on the upper half plane with finite stabilizers. Any finite order element $g\ne\pm I$ in $\slz$ has trace $0$ or $\pm1$, which is never equal to $2\mod m$ when $m\ge 4$.
\end{proof}

\begin{defn}
\label{def from torsion}
By a {multisection from torsion construction}, we mean a multisection obtained by choosing  points of type $3m$ where $m$ ranges over a nonempty finite set $I$ of positive integers. In that case, the associated cover $\xn$ is a disjoint union of $\xm$'s
\begin{equation}
    \label{torsion cover}
    \xn = \bigcup_{m\in I} \xm.
\end{equation}
The degree of $\xn$ is $9\sum_{m\in I} J_2(m)$.
\end{defn}

The multisections from torsion constructions described above are exactly those appeared in Theorem \ref{old AG}.

\section{Proofs of Theorem \ref{connected new multi} and \ref{topo constr}}
In this section, we will show that multisections of high degree are flexible: Any multisection from torsion construction of degree high enough can be deformed in two ways. More importantly, such deformations produce new multisections that are not homotopic to any ones from torsion constructions as in the previous section.  

\subsection{Two general results about deforming multisections}
\begin{prop} 
\label{deform a}
Suppose that a multisection $\sigma$ of degree $n$ is associated to a connected cover $\xn$ such that $\gn = \pi_1(\xn)/Z(K)$ is a free group. Then for any positive integer $k$, there exists another multisection $\tau$ of degree $kn$ such that 
 \begin{enumerate}
     \item the cover associated to $\tau$ has $k$ components, each isomorphic to $\xn$ as a cover of $\x$, and that
     \item as virtual sections, $\tau$ is homotopic to $k$ disjoint copies of $\sigma$ in the sense of Definition \ref{def virtual}.
 \end{enumerate}
\end{prop}


\begin{proof}
Since $\xn$ is given by a multisection, each point in $\xn$ can be uniquely represented as a pair $(F,x)$ where $x\in C_F$. 

\begin{lem}\label{nonvsec}
Let $\xn$ be as in Proposition \ref{deform a}. Let $\xi$ denote the (real) rank-2 vector bundle over $\xn$ whose fiber at $(F,x)$ is $T_{x} C_F$, the tangent space of $C_F$ at $x$. Then $\xi$ is a trivial vector bundle. 
\end{lem}
\begin{proof}
The classifying map $\phi_\xi: \xn\to B\slr$ of the vector bundle $\xi$ is a composition of two maps
$$\xn\xrightarrow[]{} B\gn \xrightarrow[]{} B\slr$$
by the same argument as in the proof of Lemma \ref{BSL}. An orientable rank-2 vector bundle is trivial if and only if its Euler class is zero. The Euler class of $\xi$ is the pullback of the universal Euler class from $B\slr$ to $\xn$ by the composition of maps above. Notice that $H^2(B\gn;\ZZ)=0$ since $\gn$ is a free group. Hence, $\xi$ has zero Euler class and therefore must be trivial. 
\end{proof} 

Hence, $\xi$ admits a smooth nonvanishing vector field $v$:
\begin{equation}
    \label{vector field v}
        \xymatrix{
T_{x} C_F \ar[r] & E_\xi \ar[d]\\
& \xn \ar@/^-10pt/[u]_{v}\\
}
\end{equation}
Equip each cubic curve $C_F$ with its unique flat Riemannian metric compatible with its complex structure with unit area. Normalize the vector field $v$ so that it has length 1 under this metric. Let $\epsilon : \x \to \RR_{>0}$ be a continuous function with the following properties: for any $F\in \x$, 
 \begin{itemize}
     \item for any $x\in\sigma(F)$, the exponential map $\exp_x: T_xC_F\to C_F$ is injective within radius $\le \epsilon(F)$, and
     \item  $\epsilon(F)<\min\{d(x,y) : \forall x,y\in\sigma(F)\text{ and } x\ne y\}$ where $d(-,-)$ is the unique flat Riemannian metric with unit volume on $C_F$.
 \end{itemize}
Such a continuous function $\epsilon$ exists because the Riemannian metric on $C_F$ varies continuously with $F\in\x$. 

For each $F\in\x$, we define
$$\tau_1(F):= \bigg\{\exp_x\Big(\frac{1}{4k}\epsilon(F) v(F,x)\Big) : x\in\sigma(F)\bigg\}.$$
The properties of $\epsilon$ and $v$ above guarantee that $\tau_1(F)$ is a set of $n$ distinct points in $C_F$ for all $F\in\x$. Furthermore, the map $\tau_1: F\mapsto \tau_1(F)$ defines a multisection of degree $n$ because the function $\epsilon$ and the vector field $v$ are all continuous. 

More generally, for any $j=1, \cdots, k$, define
\begin{align*}
    \tau_j(F)&:= \bigg\{\exp_x\Big(\frac{j}{4k}\epsilon(F) v(F,x)\Big) : x\in\sigma(F)\bigg\}\\
    \tau(F)&:=\tau_1(F)\cup\cdots\cup\tau_k(F)
\end{align*}
The properties of $\epsilon$ and $v$ above guarantee that the sets $\tau_j(F)$'s are all disjoint. Hence, $\tau(F)$ is a set of $nk$ distinct points in $C_F$. The function $\tau: F\mapsto \tau(F)$ defines a multisection of degree $nk$. See Figure \ref{vector field} for pictures of the constructions described above.
\begin{figure}[ht]
    \centering
    \includegraphics[scale=.15]{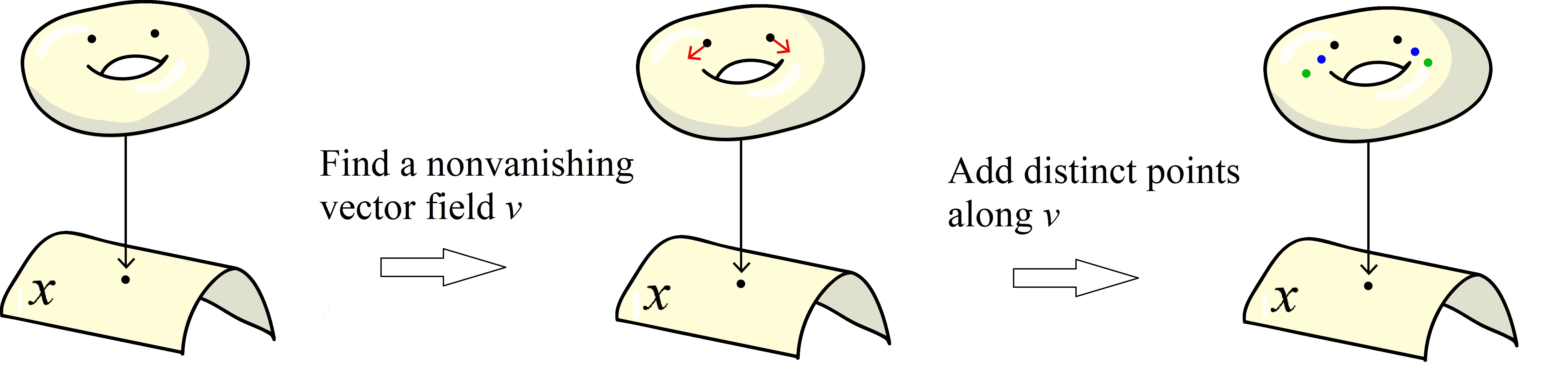}
    \caption{The figure illustrates the construction of $\tau$ when $n=2$ and $k=3$. The first picture shows a single fiber of the torus bundle $\pi:\xp\to\x$ marked with a multisection of degree $n=2$. The second picture shows a nonvanishing vector field $v$ of $\xi$. The last picture shows the new multisection $\tau$ of degree 6. Here we choose $n=2$  to simplify the  pictures. In reality, $n$ must be a multiple of 9 by Theorem \ref{old thm}.}
    \label{vector field}
\end{figure}

Finally, we check that $\tau$ satisfies the two properties stated in the proposition. Each $\tau_j$ is homotopic to $\sigma$ as multisections in the sense of Definition \ref{multi def} via the following homotopy
$$\sigma_t(F)= \bigg\{\exp_x\Big(t\frac{j}{4k}\epsilon(F) v(F,x)\Big) : x\in\sigma(F)\bigg\},\ \ \ \ \ \ t\in[0,1]$$
where $\sigma_0=\sigma$ and $\sigma_1=\tau_j$. 
Hence, $\tau = \bigcup_{j=1}^k \tau_j$ is homotopic to $k$ disjoint copies of $\sigma$ as virtual sections by Proposition \ref{mult vs virtual homotopy1}.
\end{proof}

\begin{prop} 
\label{deform b}
Suppose that a multisection $\sigma$ of degree $n$ gives a connected cover $\xn$ such that $\gn = \pi_1(\xn)/Z(K)$ is free. Then there exists a multisection $\mu$ of degree $2n$ such that
 \begin{enumerate}
     \item the cover $Y$ of $\x$ associated to the multisection $\mu$ is also a \emph{connected} 2-sheeted cover $p:Y\to \xn$, and 
     \item  $\mu$ is homotopic to $\sigma\circ p$ as virtual sections:
     $$\xymatrix{
Y \ar[r]^p \ar@/^-10pt/[rr]_{\mu \ \simeq\ \sigma\circ p} & \xn \ar[r]^\sigma & \xp.\\
}$$
    \end{enumerate}
\end{prop}
\begin{proof}
Let $\xi$ again be as in Lemma \ref{nonvsec}. We need the following lemma:
\begin{lem}\label{ntrline}
  The trivial rank-2 vector bundle $\xi$ contains a nontrivial rank-1 vector sub-bundle  over $\xn$.
\end{lem}
\begin{proof}
By Lemma \ref{nonvsec}, the total space $E$ of $\xi$ is homeomorphic to $\RR^2\times \xn$. Hence, there is a 1-1 correspondence between isomorphism classes of rank-1 sub-bundles of $\xi$ and homotopy classes of maps $f:\xn \to \RR\PP^1$. The condition that  the line bundle is trivial is equivalent to the existence of another function $g$ making the following diagram commute:
$$ \xymatrix{
& S^1\ar[d]^{2:1}\\
\xn \ar[ru]^g \ar[r]^f & \RR\PP^1\\
}$$
Suppose $ H^1(\RR\PP^1;\ZZ)= \langle \alpha\rangle$ and $H^1(S^1;\ZZ) = \langle \beta\rangle$. Then the 2-sheeted cover $S^1\to \RR\PP^1$ will take $\alpha\mapsto \pm 2\beta$. If such a lift $g$ exists, then
$$f^*(\alpha) = g^*(\beta)= \pm 2g^*(\alpha) \in 2H^1(\xn ; \ZZ).$$ 
However we can choose an $f$ such that $f^*(\alpha) \not \in 2 H^1(\xn; \ZZ)$ in the following way. Pick any element $\delta\in H^1(\xn; \ZZ)$ that is not in $2H^1(\xn; \ZZ)$. We know that such a $\delta$ exists because $H^1(\xn; \ZZ)\cong H^1(\gn;\ZZ)$  by Lemma \ref{H1 iso} where $\gn$ is a free group by hypothesis. Since $\RR\PP^1\simeq S^1$ is a $K(\ZZ, 1)$, there exists an $f$ such that $f^*(\alpha) = \delta$. This $f$ does not lift to a map to $S^1$ and hence defines a nontrivial rank-1 sub-bundle of $\xi$.
 \end{proof}
 
The rank-1 sub-bundle of $\xi$ defines a line field $\pm \omega: \xn\to E$ such that $\pm \omega(F,x)$ is a unit vector in $T_x C_F$ well-defined up to a sign for any $(F,x)\in \xn$. We will use $\pm \omega$ to construct a multisection $\mu$ of degree $2n$ (in a similar way as how we used $v$ to construct a multisection $\tau$ in the proof of Proposition \ref{deform a}). Take the same function $\epsilon: \x\to \RR_{>0}$ as in the proof of Proposition \ref{deform a}. Define 
$$\mu(F)= \bigg\{\exp_x\Big(\frac{1}{4}\epsilon(F) \omega(F,x)\Big),\  \exp_x\Big(-\frac{1}{4}\epsilon(F) \omega(F,x)\Big) : x\in\sigma(F)\bigg\}.$$
As in the proof of Proposition \ref{deform a}, properties of $\epsilon$ and $\omega$ guaranteed that $\mu(F)$ is a set of $2n$ distinct points in $C_F$. See Figure \ref{line f} for an sketch of the construction of $\mu$.
\begin{figure}[ht]
    \centering
    \includegraphics[scale=.15]{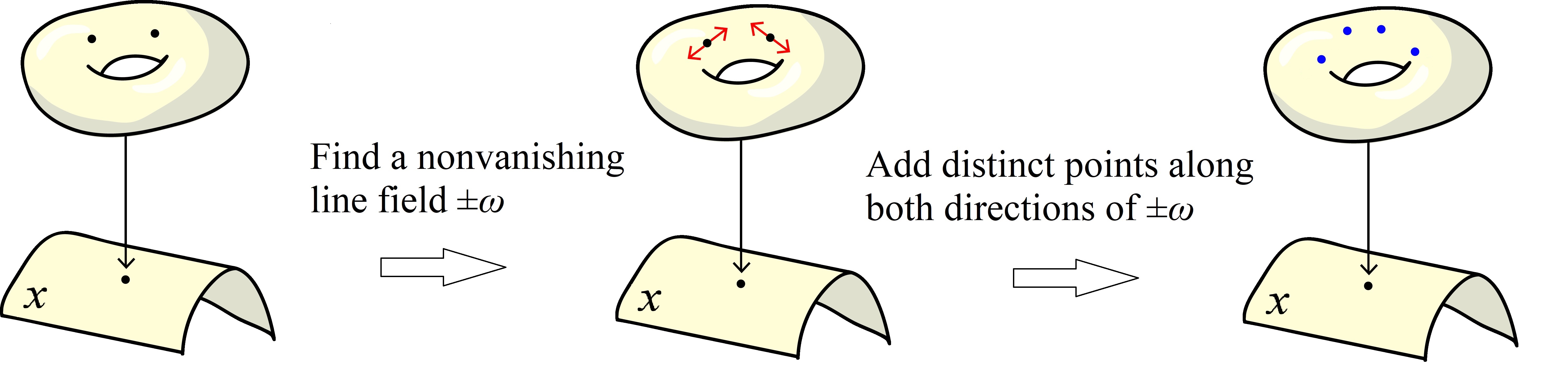}
    \caption{The figure shows the construction of $\mu$ when $n=2$. The multisection of degree $2$ in the first picture is deformed to a multisection of degree $4$ in the last picture. Here again we choose $n=2$ to simplify the pictures. In reality, $n$ must be a multiple of 9 by Theorem \ref{old thm}.}
    \label{line f}
\end{figure}

Now we check that $\mu$ satisfies the properties stated in the theorem. The cover $Y$ associated to $\mu$ is a 2-sheeted cover of $\xn$ by sending the two points 
$$\exp_x\Big(\frac{1}{4}\epsilon(F) \omega(F,x)\Big),\  \exp_x\Big(-\frac{1}{4}\epsilon(F) \omega(F,x)\Big)$$
to $x$. We first show that $Y$ is connected. Since $\xn$ is connected by hypothesis, if $Y$ is disconnected then $Y$ has two components each homeomorphic to $\xn$. Then each component will pick a sign for $\pm \omega$ at each point $(F,x)\in\xn$ and hence gives a nonvanishing vector field of the nontrivial rank-1 vector sub-bundle as in Lemma \ref{ntrline}, reaching a contradiction. 

Finally, $\mu$ is homotopic to $\sigma\circ p$ via the following homotopy:
$$\mu_t(F):= \bigg\{\exp_x\Big(\frac{t}{4}\epsilon(F) \omega(F,x)\Big),\  \exp_x\Big(-\frac{t}{4}\epsilon(F) \omega(F,x)\Big) : x\in\sigma(F)\bigg\}.$$
We have $\mu_0=\mu$ and $\mu_1=\sigma\circ p$. 
\end{proof}


Now we will prove Theorem \ref{connected new multi} and Theorem \ref{topo constr} by applying Proposition \ref{deform b} and Proposition \ref{deform a}, respectively, to deform $\sigma_m$ when $m\ge 4$ and obtain new multisections. 

\subsection{Proof of Theorem \ref{connected new multi}}
Let $\sigma_m$ be as in (\ref{sigma_m}) for $m\ge4$. By Lemma \ref{conjugate to congruence subgroup}, its associated cover $\xm$ is connected. By Corollary \ref{cong is free}, $\pi_1(\xm)/Z(K)$ is free because $m\ge4$. Hence, $\sigma_m$ satisfies the hypothesis of Proposition \ref{deform b}, which now gives us a new multisection $\mu$ of degree $2\deg(\sigma_m)=18J_2(m)$. It remains for us to prove that $\mu$ is not homotopic to any multisection from torsion construction. 



For the sake of contradiction, suppose that $\mu$ is homotopic to a multisection from torsion construction. Let $Y$ be the cover of $\x$ associated to $\mu$. By Proposition \ref{mult vs virtual homotopy1}, $Y$ is isomorphic as a cover to a disjoint union $\bigcup_{r\in I} \xn_r$ as in (\ref{torsion cover}). Since $Y$ is connected by Proposition \ref{deform b} (1), it must be the case that $Y\cong \xn_r$ for some integer $r$. Moreover, since $Y\cong \xn_r$ is a 2-sheeted cover of $\xm$ by Proposition \ref{deform b}, $\pi_1(\xn_r)$ is an index-2 subgroup of $\pi_1(\xm)$ up to conjugation in the ambient group $\pi_1(\x)$.

Now consider the images of $\pi_1(\xn_r)$ and of $\pi_1(\xm)$ under the monodromy representation $\rho:\pi_1(\x)\to\slz$ as in (\ref{monodromy rep}). By the discussion above,  $\rho(\pi_1(\xn_r))$ is an index-2 subgroup of $\rho(\pi_1(\xm))$ up to conjugation in $\slz$. By Lemma \ref{conjugate to congruence subgroup}, we know that $\rho(\pi_1(\xn_r))$ is conjugate to $\Gamma_1(r)$ in $\slz$, while $\rho(\pi_1(\xm))$ is conjugate to $\Gamma_1(m)$. Now we have that $\Gamma_1(r)$ is conjugate to an index-2 subgroup of $\Gamma_1(m)$ in $\slz$, which implies that $m|r$ by comparing the traces of matrices in $\Gamma_1(r)$ and in $\Gamma_1(m)$. Moreover, by comparing the degrees of $\mu$ and of $\sigma_m$, we have
\begin{equation}
    \label{local1}
    \deg(\mu) = 2\deg(\sigma_m) = 18 J_2(m) = \deg(\sigma_r) = 9J_2(r) \ \ \Longrightarrow\ \  2J_2(m)=J_2(r).
\end{equation}
Observe that the formula of Jordan's 2-totient function $J_2$ in (\ref{number of points of type 3m}) implies that
\begin{equation}
    \label{local2}
    J_2(r)\ge J_2(m)\cdot J_2(r/m).
\end{equation}
(\ref{local1}) implies that $r/m>1$. Moreover, (\ref{local1}) and (\ref{local2}) together imply that $2\ge J_2(r/m)$, which is impossible by the formula of $J_2$ in (\ref{number of points of type 3m}). 
\qed

\subsection{Proof of Theorem \ref{topo constr}}
Suppose that $n$ is as in Theorem \ref{topo constr}:
$$n=9\sum_{m\in I} k_m J_2(m)$$
where $I$ is a finite set of positive integers and each $k_m$ is a positive integer such that $k_m=1$ for every $m\le3$ and that $k_m>1$ for some $m\ge 4.$ We need to produce a multisection $\tau$ of degree $n$ and prove that $\tau$ is not homotopic to any multisection from torsion construction. 

We will proceed by modifying the proof of Proposition \ref{deform a}. Let $\epsilon:\x\to\RR_{>0}$ be a continuous function such that for any $F\in \x$,
 \begin{itemize}
     \item  for any $m\in I$ and for any $x\in\sigma_m(F)$, the exponential map $\exp_x: T_xC_F\to C_F$ is injective within radius $\le \epsilon(F)$, and
     \item  $\epsilon(F)<\min\{d(x,y) : x,y\in\bigcup_{m\in S}\sigma_m(F)\text{ and } x\ne y\}$ where $d(-,-)$ is the unique flat Riemannian metric on $C_F$ with unit volume.
 \end{itemize}

For each $m\in I$ such that $k_m>1$, it is necessary that $m\ge 4$ by hypothesis. Hence, $\sigma_m$ satisfies the assumption of Proposition \ref{deform a} by the same argument as in the first paragraph of the proof of Theorem \ref{connected new multi} above. Similar to the proof of Proposition \ref{deform a}, we define $v_m$ to be a smooth unit vector field of the following trivial vector bundle:
$$    \xymatrix{
T_{x} C_F \ar[r] & E \ar[d]\\
& \xm \ar@/^-10pt/[u]_{v_m}\\
}$$
In fact, we can choose $v_m$ to be the composition of $v$ as in (\ref{vector field v}) and the covering map $\xm\to\x$. 
For each $m\in I$ such that $k_m>1$, we define
$$    \tau_m(F)= \bigg\{\exp_x\Big(\frac{j}{4k}\epsilon(F) v(F,x)\Big) : x\in\sigma_m(F),\ j=1,...,k_m\bigg\}.
$$ 
For each $m\in I$ such that $k_m>1$, we define $\tau_m:=\sigma_m$. Finally, define
$$\tau(F):=\bigcup_{m\in I} \tau_m(F).$$
By the properties of $\epsilon$ and $v_m$, the set $\tau(F)$ contains $n$ distinct points for every $F\in\x$. The function $\tau:\x\to\xu$ is a multisection of degree $n$.

It remains to show that $\tau$ is not homotopic to any multisection from torsion construction. By assumption, there exists some $m\in I$ such that $k_m>1$. It suffices to prove that $\tau_m$ is not homotopic to any multisection $\sigma$ from torsion construction. Suppose not. By Proposition \ref{deform a}, the cover associated to $\tau_m$ has $k$ components, each homeomorphic to $\xm$. By (\ref{torsion cover}), the cover associated to $\sigma$ has the following form
\begin{equation}
    \label{S}
    \bigcup_{r\in S} \xn_{r}
\end{equation}
where $S$ is a finite set of positive integers. By Proposition \ref{mult vs virtual homotopy1}, since $\sigma$ is homotopic to $\tau_m$ as multisections, each component $\xn_r$ in (\ref{S}) must be isomorphic to $\xm$ as covers of $\x$. From this, we claim that it must be that $r=m$. Notice that it is possible for $\xn_r$ and $\xm$ to have the same degree yet $r\ne m$ (for example, $\xn_{5}$ and $\xn_{6}$ both have degree 216), so we need some finer invariants than degrees. To achieve our goal, we again consider the images of $\pi_1(\xn_r)$ and of $\pi_1(\xm)$ under the monodromy representation $\rho:\pi_1(\x)\to\slz$. Since $\xn_r\cong \xm$ as covers of $\x$, the two subgroups $\pi_1(\xn_r)$ and $\pi_1(\xm)$ are conjugate in $\pi_1(\x)$. Moreover, by Lemma \ref{conjugate to congruence subgroup}, we know that $\rho(\pi_1(\xn_r))$ is conjugate to $\Gamma_1(r)$ in $\slz$, while $\rho(\pi_1(\xm))$ is conjugate to $\Gamma_1(m)$. Since  $\Gamma_1(r)$ and $\Gamma_1(m)$ are conjugate in $\slz$, we conclude that $r=m$ by comparing the traces of matrices in $\Gamma_1(r)$ and $\Gamma_1(m)$. Hence, the indexing set $S$ is just a singleton $\{m\}$ and $\sigma=\sigma_m$. However, this is in contradiction with our initial assumptions that $\sigma$ is homotopic to $\tau_m$ and that $\tau_m$ is of degree $k\deg(\sigma_m)$ for $k>1$. \qed

\begin{rem}
One might wonder: What will go wrong if do the same proofs of Theorem \ref{connected new multi} and \ref{topo constr} above but for $m<4$? It turns out that $\Gamma_1(m)$ is not free when $m<4$, because its Euler characteristic 
$$\chi\Big(\Gamma_1(m)\Big)=J_2(m)\chi\Big(\slz\Big)=-\frac{J_2(m)}{12}$$
is not an integer. Moreover, one can check that when $m<4$, the vector bundle $\xi$ as in Lemma \ref{nonvsec} is nontrivial and hence our constructions above using nonvanishing vector fields would fail. Thus, to improve Theorem \ref{connected new multi} and \ref{topo constr}, one must find a different topological construction beyond what is presented in this paper.
\end{rem}

\section{Appendix}
In this Appendix, we record McMullen's observation that every algebraic multisection of $\pi$ must be a multisection from torsion construction (Theorem \ref{algebraic is torsion} below). None of the results in this Appendix are due to the authors.

The space $\x$ of smooth cubic curves is a complex quasi-projective variety because it is the complement of the discriminant hypersurface in $\CC\PP^{9}$. The universal cubic curve $\pi:\xp\to\x$ is a morphism of algebraic varieties. We say a multisection $\sigma$ of $\pi$ is \emph{algebraic} if it is also a morphism of algebraic varieties. 

\begin{thm}
\label{algebraic is torsion}
Every algebraic multisection of $\pi$ must be a multisection from torsion construction.
\end{thm}

The proof will proceed in the following three steps.

\subsection{Holomorphic sections and the Teichm\"uller spaces}
Let $\TT$ denote the Teichm\"uller space of marked Riemann surfaces of genus 1 and with $1$ puncture. It is well-known that the following map gives a biholomorphism between $\TT$ and the upper half plane $\HH$:
\begin{align*}
    \HH&\longrightarrow\TT\\
    \tau&\longmapsto E_\tau:=\CC/\langle 1,\tau\rangle \ \ \ \text{ with a marked basis $\{1,\tau\}$}
\end{align*}
The \emph{universal punctured elliptic curve} $\CT\to\TT$ is a bundle whose fiber over a point $\tau\in \TT$ is the \underline{punctured} Riemann surface $E_\tau\setminus\{0\}$. 

\begin{lem}\label{royden}
Every holomorphic section of the universal punctured elliptic curve $\CT\to \TT$ is of the form $s(\tau)=a+b\tau$ for a fixed $(a,b)\in \RR^2\setminus \ZZ^2$ independent of $\tau$.
\end{lem}

\begin{proof}
Let $\mathcal{T}_{1,2}$ denote the Teichm\"uller  space of marked Riemann surfaces of genus 1 and with 2 punctures. There is a natural projection $p:\mathcal{T}_{1,2}\to \TT$ defined by forgetting the second marked point.  Note that $\mathcal{T}_{1,2}$ is the universal cover of $\CT$. Thus, any holomorphic section $s:\TT\to\CT$ lifts to a holomorphic section $\tilde s:\TT\to\mathcal{T}_{1,2}$. Consider the following two maps:
$$\xymatrix{
\mathcal{T}_{1,2} \ar[r]_p & \TT \ar@{-->}@/^-10pt/[l]_{\tilde s}
}$$
Both $p$ and $\tilde s$ are holomorphic and are therefore distance non-increasing with respect to the Kobayashi metrics on the two Teichm\"uller spaces. By Royden's theorem \cite{Royden}, the Teichm\"uller metric and Kobayashi metric coincide. Thus, $\tilde s$ is an isometric embedding with respect to the Teichm\"uller metrics on its domain and codomain. Now, it follows from the uniqueness of Teichm\"uller maps that $\tilde s$ must be of the form
$$\forall\tau\in\TT \ \ \ \ \ \ \ \ \ \ \tilde s(\tau)=a+b\tau$$
for some fixed $(a,b)\in \RR^2\setminus\ZZ^2$ independent of $\tau$. The section $s:\TT\to\CT$ is nothing but $\tilde s:\TT\to\mathcal{T}_{1,2}$ after forgetting the marked basis of the lattice. 
\end{proof}

\subsection{Algebraic multisections and the Legendre family}
Consider the Legendre family of cubic plane curves:
$$C_t: y^2=x(x-1)(x-t) \ \ \ \ \ \ t\in \CC \setminus \{0,1\}.$$
We assume that every cubic curve $C_t$ includes a point at infinity which is a flex point and serves as an identity for its elliptic curve group law. 
\begin{prop}
\label{Legendre}
Suppose that $\sigma$ is an algebraic section of the  bundle
$$\uconf_n C_t\to \xu\to \CC\setminus\{0,1\}.$$
Then every element of the finite set $\sigma(t)$ must be a torsion point in $(C_t,\infty)$ for every $t\in \CC\setminus\{0,1\}.$
\end{prop}
\begin{proof}
We will first prove the following claim, only assuming that $\sigma$ is holomorphic.
\begin{claim}
Suppose that $\sigma$ is a holomorphic section of the  bundle
$$\uconf_n C_t\to \xu\to \CC\setminus\{0,1\}.$$
Then the difference of any two elements in $\sigma(t)$ is a torsion point in $(C_t,\infty)$ for every $t\in \CC\setminus\{0,1\}.$
\end{claim}
\begin{proof}
Given such a $\sigma$, define the following space
$$Y:=\Big\{(t,x,y) : t\in \CC\setminus\{0,1\},\ \text{$x$ and $y$ are two distinct elements in $\sigma(t)$} \Big\}.$$
Since $\sigma$ is a multisection, $Y$ is a finite cover of $\CC\setminus\{0,1\}$ by projecting onto the first coordinate. Let $Z$ be any path component of $Y$. Since every elliptic curve is isomorphic to a cubic curve in the Legendre family, we have a natural map $\TT\to\CC\setminus\{0,1\}$ which gives the universal cover of $\CC\setminus\{0,1\}$ and of $Z$. The map 
\begin{align*}
    s: Z&\to C_t\setminus\{\infty\}\\
    (t,x,y)&\mapsto x-y \ \ \ \text{ under the group law of $(C_t,\infty)$}
\end{align*}
defines a section of the (punctured) universal curve over $Z$. This section pulls back to the universal cover $\TT$ of $Z$ and gives a holomorphic section $\tilde s:\TT\to\CT$ as in Lemma \ref{royden}. Hence, we have that $\tilde s(\tau)=a+b\tau$ for some $(a,b)\in  \RR^2$. Finally, under the natural action of $\slz$ on the marked elliptic curve $\CC/(\ZZ+\tau\ZZ)$, the stabilizer of $\tilde s(\tau)$ is exactly $\pi_1(Z)$. Since $Z$ is a connected finite cover of $\CC\setminus\{0,1\}$, we have that
$$[\pi_1(Z):\slz]\le [\pi_1(\CC\setminus\{0,1\}):\slz]<\infty.$$
We thus conclude that the $\slz$-orbit of $\tilde s(\tau)$ is finite. This implies that $(a,b)\in\QQ^2$ and that $\tilde s(\tau)$ is a torsion point of $\CC/(\ZZ+\tau\ZZ)$ for every $\tau\in \TT$. Since $\tilde s$ is a lift of $s$, we conclude that $s(t,x,y)$ is always a torsion point of $(C_t,\infty)$.
\end{proof}

By this claim, we have that every holomorphic section $\sigma$ of the bundle
$$\uconf_n C_t\to \xu\to \CC\setminus\{0,1\}$$
locally is of the following form
$$\sigma(t)=p(t)+\Sigma(t)$$
where $p(t)$ is an arbitrary point on $C_t$ and $\Sigma(t)$ is a finite set of torsion points on $(C_t,\infty)$. 
Hence, there exists an integer $N$ such that 
$N\cdot \sigma(t)$ is a single point on $C_t$ for every $t\in \CC\setminus\{0,1\}$. So the function $s(t):=N\cdot \sigma(t)$ defines a section of the universal Legendre curve 
\begin{equation}
    \label{leg}
    \xymatrix{
C_t \ar[r] & E \ar[r] & \CC\setminus\{0,1\} \ar@{-->}@<-1ex>@/^-10pt/[l]_{s}\\
}
\end{equation}

If $\sigma$ is an algebraic multisection, then the $s$ that we constructed above is an algebraic section. Algebraic sections of (\ref{leg}) form the Mordell-Weil group of an elliptic surface given by the Legendre family. It is well-known that there are only four such algebraic sections, given by  $s(t)=(x,y)=(0,0),(1,0),(t,0),$ or $\infty$. See, for example, Theorem 8.9 in Section 8.3 of \cite{SS}. In particular,  $s(t)$ is always a 2-torsion of $(C_t,\infty)$. Hence, every point in the multisection $\sigma(t)$ must also be a torsion point on $(C_t,\infty)$.
\end{proof}

\subsection{Finishing the proof of Theorem \ref{algebraic is torsion}}
\begin{proof}[Proof of Theorem \ref{algebraic is torsion}]
Suppose that $\sigma:\x\to\xu$ is an algebraic multisection. It is a standard fact that every smooth cubic curve can be mapped by a projective linear transformation to a curve in a Legendre family. Since a projective linear map takes flex points to flex points,  Proposition \ref{Legendre} implies that every point in $\sigma(F)$ must be a torsion point on $C_F$ with respect to some flex point as identity.

Finally, we check that $\sigma$ must be a multisection from torsion construction in the sense of Definition \ref{def from torsion}. Let $\xn$ denote the cover associated to the multisection $\sigma$. Take any $(F,x)\in\xn$. Let $N$ be the smallest number such that $x$ is an $N$-torsion on $C_F$ with any one of the nine flex points as identity. Since the nine flex points differ by 3-torsion translations, we must have that $N=3m$ for some integer $m$ and hence $x$ is a point of type $3m$ as in Definition \ref{type 3m}. We have proven in Lemma \ref{conjugate to congruence subgroup} that the cover $\xm$, consisting of a cubic curve and a point of type $3m$ as defined in (\ref{xm define}) above, is connected. Hence, the connected component of $\xn$ containing $(F,x)$ must be the entire $\xm$.  Hence, $\sigma$ must be a multisection from torsion construction as in Definition \ref{def from torsion}.
\end{proof}




\end{document}